\documentclass[12pt,reqno]{amsart}

\title[]{From the Schr\"odinger problem to the Monge-Kantorovich problem}
\author{Christian L\'eonard}
\date{November 2010}

\usepackage{amssymb, amsmath, amsfonts, latexsym, enumerate}

\newtheorem{theorem}{Theorem}
\newtheorem{lemma}[theorem]{Lemma}
\newtheorem{proposition}[theorem]{Proposition}
\newtheorem{corollary}[theorem]{Corollary}

\theoremstyle{remark}
\newtheorem{remark}[theorem]{Remark}
\newtheorem{remarks}[theorem]{Remarks}
\newtheorem{example}[theorem]{Example}
\newtheorem{examples}[theorem]{Examples}

\numberwithin{theorem}{section}

\newcommand{\RR}{\mathbb{R}}
\newcommand{\PP}{\mathbb{P}}

\newcommand{\1}{\textbf{1}}

\newcommand\pf{_{\#}}

    \DeclareMathOperator{\dom}{dom}
    
    \DeclareMathOperator{\cl}{cl}
    \DeclareMathOperator{\cv}{cv}
    \DeclareMathOperator{\inter}{int}
    
    \DeclareMathOperator{\supp}{supp}
    \DeclareMathOperator{\ls}{ls}
    
    \DeclareMathOperator{\argmin}{argmin}

\newcommand{\Boulette}[1]{\par\medskip\noindent $\bullet$\ Proof of #1.}

\newcommand\seq[2]{(#1_#2)_{#2\ge1}}
\newcommand\Seq[2]{(#1^#2)_{#2\ge1}}
\newcommand\Lim[1]{\lim_{#1\rightarrow\infty}}
\newcommand\Liminf[1]{\liminf_{#1\rightarrow\infty}}
\newcommand\Limsup[1]{\limsup_{#1\rightarrow\infty}}
\newcommand\Glim[1]{\Gamma\textrm{-}\lim_{#1\rightarrow\infty}}

\newcommand{\eqlaw}{\overset{\textrm{Law}}=}

\newcommand{\lsc}{lower semicontinuous}

\newcommand\pont[3]{#1^{#2,#3}}
\newcommand\haut[3]{#1^{#2#3}}
\newcommand\bas[3]{#1_{#2#3}}

\newcommand\XX{\mathcal{X}}
\newcommand\XXX{\XX^2}
\newcommand\PX{\mathrm{P}(\XX)}
\newcommand\PXX{\mathrm{P}(\XXX)}
\newcommand\PO{\mathrm{P}(\Omega)}
\newcommand\CX{\mathrm{C}_b(\XX)}

\newcommand\CXX{\mathrm{C}_b(\XXX)}
\newcommand\OO{\Omega}
\newcommand\Oac{\Omega_{\textrm{ac}}}
\newcommand\CO{\mathrm{C}_b(\Omega)}
\newcommand\MO{\mathrm{M}_b(\Omega)}
\newcommand\MOp{\mathrm{M}_b^+(\Omega)}

\newcommand\ii{{[0,1]}}
\newcommand\IX{\int_{\XX}}
\newcommand\IXX{\int_{\XXX}}
\newcommand\IO{\int_\Omega}
\newcommand\Iii{\int_\ii}
\newcommand\Ph{\widehat{P}}
\newcommand\pih{\widehat{\pi}}

\newcommand{\infX}{\inf_{x\in X}}
\newcommand{\infY}{\inf_{y\in Y}}
\newcommand{\infV}{\inf_{x\in V}}
\newcommand{\supY}{\sup_{y\in Y}}

\newcommand{\xy}{\langle x,y\rangle}

\begin{document}

% ************** page de garde ******************************

 \address{Modal-X. Universit\'e Paris Ouest. B\^at.\! G, 200 av. de la R\'epublique. 92001 Nanterre, France}
 \email{christian.leonard@u-paris10.fr}

\begin{abstract}
The aim of this article is to show that the Monge-Kantorovich
problem is the limit of a sequence of entropy minimization
problems when a fluctuation parameter tends down to zero. We prove
the convergence of the entropic values to the optimal  transport
cost as the fluctuations decrease to zero, and we  also show that
the limit points of the entropic minimizers are optimal transport
plans. We investigate the dynamic versions of these problems by
considering random paths and describe the connections between the
dynamic and static problems. The proofs are essentially based on
convex and functional analysis. We also need specific properties
of $\Gamma$-convergence which we didn't find in the literature.
Hence we  prove these $\Gamma$-convergence results which are
interesting in their own right.
\end{abstract}

\maketitle \tableofcontents

% ************* corps du texte ****************************

\section{Introduction}

The aim of this article is to describe a link between the
Monge-Kantorovich optimal transport problem and a sequence of
entropy minimization problems. We show that the Monge-Kantorovich
problem is the limit of this sequence when a fluctuation parameter
tends down to zero. More precisely, we prove that the entropic
values tend to the optimal cost as the fluctuations decrease to
zero, and also that the limit points of the entropic minimizers
are optimal transport plans. We also investigate the dynamic
versions of these problems by considering random paths.

Our main  results are  stated at Section \ref{sec-results}, they
are Theorems \ref{res-07}, \ref{res-04} and \ref{res-09}.

Although the assumptions of these results are in terms of large
deviation principle, it is not necessary to be acquainted to this
theory or even to probability theory to read this article. It is
written for analysts and we tried as much as possible to formulate
the probabilistic notions in terms of analysis and measure theory.
A short reminder of the basic definitions and results of large
deviation theory is given at the Appendix.

In its Kantorovich form, the optimal transport problem dates back
to the 40's, see \cite{Kanto42, Kanto48}. It appears that its
entropic approximation has its roots in an even older problem
which was addressed by Schr\"odinger in the early 30's in
connection with the newly born wave mechanics, see \cite{Sch32}.

The Monge-Kantorovich optimal transport problem is about finding
the best way of transporting some given mass distribution  onto
another one.  We describe these mass distributions by means of two
probability measures on a state space $\XX:$ the initial one is
called $\mu_0\in\PX$ and the final one $\mu_1\in\PX$ where $\PX$
is the set of all probability measures on $\XX.$ The rules of the
game are (i): it costs $c(x,y)\in[0,\infty]$ to transport a unit
mass from $x$ to $y$ and (ii): it is possible to transport
infinitesimal portions of mass from the $x$-configuration $\mu_0$
to the $y$-configuration $\mu_1.$ The resulting minimization
problem is the celebrated Monge-Kantorovich problem
\begin{equation}\label{eq-36}
     \IXX c\,d\pi \rightarrow \min;\quad \pi\in\PXX : \pi_0=\mu_0,
    \pi_1=\mu_1
\end{equation}
where $\PXX$ is the set of all probability measures on $\XXX$ and
$\pi_0,\pi_1\in\PX$ are respectively the first and second marginal
measures of the joint probability measure $\pi\in\PXX.$ Optimal
transport is an active field of research. For a remarkable
overview of this exciting topic, see Villani's textbook
\cite{Vill09} and the references therein.

Now, let us have a look at Schr\"odinger's problem. Suppose that
you observe a very large number of non-interacting
indistinguishable particles which are approximately distributed
around a probability measure $\mu_0\in\PX$ on the state space
$\XX.$ We view $\mu_0$ as the initial configuration of the whole
particle system. Suppose that you know that the dynamics of each
individual particle is driven by a stochastic process whose law is
$R^k\in\PO:$ i.e.\ a probability measure on the space
$$\OO=\XX^{\ii}$$ of all paths\footnote{During the rigorous
treatment, we shall only consider subspaces $\OO$ of $\XX^{\ii},$
for instance the subspace of all continuous paths.} from the time
interval $\ii$ to the state space $\XX.$ The parameter $k$
describes the fluctuation level $1/k.$ As $k$ tends to infinity,
$R^k$ tends to some deterministic dynamics: $R^\infty$ describes a
deterministic flow. As a typical example, one can think of $R^k$
as the law of a Brownian motion with diffusion coefficient $1/k.$
Knowing this dynamics, the law of large numbers tells you that you
should expect to see the configuration of the large particle
system at the final time $t=1$ not very far from some expected
configuration, with a very high probability. Now, suppose that you
observe the system in a configuration close to some $\mu_1\in\PX$
which is far from the expected one. Schr\"odinger's question is :
\emph{``Conditionally on this very rare event, what is the most
likely path of the whole system between the times $t=0$ and
$t=1$?"} As will be seen at Section \ref{sec-problems}, the answer
to this question is related to the entropy minimization problem
\begin{equation}\label{eq-35}
   \frac 1k H(P|R^k)\rightarrow \min;\quad  P\in\PO: P_0=\mu_0, P_1=\mu_1
\end{equation}
where $ H(P|R^k)$ is the relative entropy of $P$ with respect to
the reference stochastic process $R^k$ and the renormalization
factor $1/k$ is here to prevent the entropy from exploding as the
fluctuations of $R^k$ decrease. Recall that
    $
H(P|R):=\IO\log\left(\frac{dP}{dR}\right)\,dP \in[0,\infty],
    $
for any $P,R\in\PO.$ Schr\"odinger's problem looks like the
Monge-Kantorovich one not only because of $\mu_0$ and $\mu_1,$ but
also because of some cost of transportation. Indeed, if the random
dynamics creates a trend to move in some direction rather than in
another one, it costs less to the particle system to end up at
some configurations $\mu_1$ than others. Even if no direction is
favoured, we shall see that the structure of the random
fluctuations which is described by the sequence $\Seq Rk$ encodes
some zero-fluctuation cost function $c$ on $\XXX.$

Remark that, although $1/k$ should be of the order of Planck's
constant $\hbar$ to build a Euclidean analogue of the quantum
dynamics, in \cite{Sch32} Schr\"odinger isn't concerned with the
semiclassical limit $k\rightarrow \infty.$ Let us also mention
that Schr\"odinger's paper is the starting point of the Euclidean
quantum mechanics which was developed by Zambrini \cite{CZ08}.

\subsection*{An informal presentation of the convergence result}
Our assumption is that $\Seq Rk$ satisfies a large deviation
principle in
 the path space $\OO.$ This roughly means that
\begin{equation}\label{eq-37}
    R^k(A)\underset{k\rightarrow\infty}\asymp \exp\left(-k\inf_{\omega\in
    A}C(\omega)\right),
\end{equation}
for some rate function $C:\Omega\to[0,\infty]$ and a large class
of measurable subsets $A\in\OO.$ For a rigorous definition of a
large deviation principle and basic results about large deviation
theory, see the Appendix. Very informally, the most likely paths
$\omega$ correspond to high values $R^k(d\omega)$ and therefore to
low values of $C(\omega).$ Under endpoint constraints, it
shouldn't be surprising to meet the following family of geodesic
problems
\begin{equation*}
  C(\omega)\rightarrow \min;\quad  \omega\in\OO: \omega_0=x,
  \omega_1=y
\end{equation*}
where $x,y$ describe $\XX$ and $\omega_0$ and $\omega_1$ are the
initial and final positions of the path $\omega.$ We see that the
large deviation behavior of the sequence $\Seq Rk$ brings us a
family of geodesic paths. It will be shown that the limit (in some
sense to be made precise) of the problems \eqref{eq-35} is the
Monge-Kantorovich problem with the ``geodesic" cost function
\begin{equation}\label{eq-38}
     c(x,y)=\inf\{C(\omega);\omega\in\Omega:
    \omega_0=x,\omega_1=y\},
    \quad x,y\in\XX.
\end{equation}
For instance, if $R^k$ is the law of a Brownian motion on
$\XX=\RR^d$ with diffusion coefficient $1/k,$ the rate function
$C$ is given by Schilder's theorem, a standard large deviation
result which tells us that $C$ is the classical kinetic action
functional which is given for any path $\omega$ by
$$
    C(\omega)=\frac 12 \Iii |\dot\omega_t|^2\,dt\in[0,\infty]
$$
if $\omega=(\omega_t)_{0\le t\le1}$ is absolutely continuous
($\dot\omega$ is its time derivative), and $+\infty$ otherwise.
The corresponding static cost is the standard quadratic cost
$$c(x,y)=\frac 12|y-x|^2,\quad x,y\in\RR^d.$$
As a consequence of our general results, we obtain that if the
quadratic cost transport problem admits a unique solution, for
instance when $\IX |x|^2\mu_0(dx), \IX |y|^2\mu_1(dy)<\infty$ and
$\mu_0$ is absolutely continuous, then the sequence $\Seq\Ph k$
built with the diffusion processes which are the unique solutions
to \eqref{eq-35} as $k$ varies, converges to the deterministic
process
$$
    \Ph(\cdot)=\IXX \delta_{\haut\sigma xy}(\cdot)\,\pih(dxdy)\in\PO
$$
where for each $x,y\in\XX,$ $\haut\sigma xy$ is the constant
velocity geodesic between $x$ and $y,$ $\delta_{\haut\sigma xy}$
is the Dirac measure at $\haut\sigma xy$ and $\pih\in\PXX$ is the
unique solution to the quadratic cost Monge-Kantorovich transport
problem \eqref{eq-36}. The marginal flow of $\Ph$ is defined to be
$(\Ph_t)_{0\le t\le 1}$ where for each $0\le t\le1,$ $\Ph_t=\IXX
\delta_{\haut\sigma xy_t}(\cdot)\,\pih(dxdy)\in\PX$ is the law of
the random position at time $t$ when the law of the whole random
path is $\Ph\in\PO.$ This flow is precisely the \emph{displacement
interpolation} between $\mu_0$ and $\mu_1$ with respect to the
quadratic cost transport problem, see \cite[Chapter 7]{Vill09} for
this notion.

\subsection*{Presentation of the results}

The quadratic cost is an important instance of transport cost, but
our results are valid for any cost functions $c$ and $C$
satisfying \eqref{eq-37} and \eqref{eq-38}, plus some coercivity
properties. For each $k\ge1,$ denote $\rho^k\in\PXX$ the law of
the couple of initial and final positions of the random path
driven by $R^k\in\PO.$ Then,
\begin{equation}\label{eq-39}
     \frac1k H(\pi|\rho^k) \rightarrow \min;\quad \pi\in\PXX : \pi_0=\mu_0,
    \pi_1=\mu_1
\end{equation}
is the static ``projection'' of \eqref{eq-35}.

In the sequel, any limit of sequences of probability measures is
understood with respect to the usual narrow topology. Theorem
\ref{res-09} states that, as $k$ tends to infinity,  there exists
a sequence $\Seq{\mu_1}k$ in $\PX$ such that $\Lim k\mu_1^k=\mu_1$
and the modified minimization problem
\begin{equation}\label{eq-39b}
     \frac1k H(\pi|\rho^k) \rightarrow \min;\quad \pi\in\PXX : \pi_0=\mu_0,
    \pi_1=\mu_1^k
\end{equation}
verifies the following two assertions:
\begin{itemize}
    \item The minimal value of \eqref{eq-39b} tends to the minimal
    value of \eqref{eq-36}, where $c$ is given by \eqref{eq-38},
    which is precisely the optimal transport cost $T_c(\mu_0,\mu_1)$;
    \item If $T_c(\mu_0,\mu_1)$ is finite, for all large enough
    $k,$ \eqref{eq-39b} admits a unique minimizer $\pih^k,$ the
    sequence $\Seq\pih k$ admits limit points in $\PXX$ and any such limit point
    is a solution to the Monge-Kantorovich problem \eqref{eq-36},
    i.e.\ an optimal transport plan.
\end{itemize}
It is not necessary that $c$ is derived from a dynamical cost $C$
via \eqref{eq-38}. A similar result holds in this more general
setting, this is the content of Theorem \ref{res-07}. The
dynamical analogue of this convergence result is stated at Theorem
\ref{res-04} and the connection between the dynamic and static
minimizers is described at Theorem \ref{res-09}.

Examples of random dynamics $\Seq Rk$ are introduced. They are
mainly based on random walks so that one can compute the
corresponding cost functions $C$ and $c.$ In particular, we
propose dynamics which generate the standard costs
$c_p(x,y):=|y-x|^p,$ $x,y\in\RR^d$ for any $p>0,$ see Examples
\ref{ex-02} for such dynamics based on the Brownian motion.

We also prove technical results about $\Gamma$-convergence which
we didn't find in the literature. They are efficient tools for the
proofs of the above mentioned convergence results. A typical
result about the $\Gamma$-convergence of a sequence of convex
functions $\seq fk$ is: If the sequence of the convex conjugates
$\seq{f^*}k$ converges pointwise, then $\seq fk$
$\Gamma$-converges. Known results of this type are usually stated
in separable reflexive Banach spaces, which is a natural setting
when working with PDEs. But here, we need to work with the narrow
topology on the set of probability measures. Theorem \ref{res-41}
is such a result in this weak topology setting.

Finally, we also proved Theorem \ref{res-10} which tells us that
if one adds a continuous constraint to an equi-coercive sequence
of $\Gamma$-converging minimization problems, then the minimal
values and the minimizers of the new problems still enjoy nice
convergence properties.

\subsection*{Literature} The connection between  large
deviation and optimal transport has already been done by Mikami
\cite{Mikami04} in the context of the quadratic transport.
Although no relative entropy  appears in \cite{Mikami04} where an
optimal control approach is performed, our results might be seen
as extensions of Mikami's ones. In the same spirit, still using
optimal control,  Mikami and Thieullen \cite{MT06,MT08} obtained
Kantorovich type duality results.

Recently, Adams, Dirr, Peletier and Zimmer \cite{ADPZ10} have
shown that the small time large deviation behavior of a large
particle system is equivalent up to the second order to a single
step of the Jordan-Kinderleher-Otto gradient flow algorithm. This
is reminiscent of the Schr\"odinger problem, but the connection is
not completely understood by now.

The connection between the Monge-Kantorovich and the Schr\"odinger
problems is also exploited implicitly in some works where
\eqref{eq-36} is penalized by a relative entropy, leading to the
minimization problem
\begin{equation*}
      \IXX c\,d\pi +\frac1k H(\pi|\rho) \rightarrow \min;\quad \pi\in\PXX : \pi_0=\mu_0,
    \pi_1=\mu_1
\end{equation*}
where $\rho\in\PXX$ is a fixed reference probability measure on
$\XXX,$ for instance $\rho=\mu_0\otimes\mu_1.$ Putting
$\rho^k=Z_k^{-1}e^{-kc}\,\rho$ with $Z_k=\IXX
e^{-kc}\,d\rho<\infty,$ up to the additive constant $\log(Z_k)/k,$
this minimization problem rewrites as \eqref{eq-39}. See for
instance the papers by R\"uschendorf and Thomsen \cite{RT93,RT98}
and the references therein. Also interesting is the recent paper
by Galichon and Salanie \cite{GS10} with an applied point of view.

Proposition \ref{res-02} below is an important technical step on
the way to our main results. A variant of this proposition, under
more restrictive assumptions than ours, was proved by Dawson and
G\"artner \cite[Thm 2.9]{DG94} in a context which is different
from optimal transport and with no motivation in this direction.
Indeed, \cite{DG94} is aimed at studying the large deviations of a
large number of diffusion processes subject to a hierarchy of
mean-field interactions, by means of random variables which live
in $\mathrm{P}(\PO)$: the set of probability measures on the set
of probability measures on the path space $\OO$. The proofs of
Proposition \ref{res-02} in the present paper and in \cite{DG94}
differ significantly. Dawson-G\"artner's proof is essentially
probabilistic while the author's one is analytic. The strategy of
the proofs are also separate: Dawson-G\"artner's proof is based on
rather precise probability estimates which partly rely on the
specific structure of the problem, while the present one takes
place in the other side of convex duality, using the
Laplace-Varadhan principle and $\Gamma$-convergence. Because of
these significantly different proofs and of the weakening of the
hypotheses in the present paper, we provide a complete analytic
proof of Proposition \ref{res-02} at Section \ref{sec-proofs}.

\subsection*{Organization of the paper}

Section \ref{sec-problems} is devoted to the presentations of the
Monge-Kantorovich and Schr\"odinger problems. We also show
informally that they are tightly connected. Our main results are
stated at Section \ref{sec-results}. We also give here a simple
illustration of these abstract results by means of Schr\"odinger's
original example based on the Brownian motion. Since our primary
object is the sequence of random processes $\Seq Rk$, it is
necessary to connect it with the cost functions $C$ and $c.$ This
is the purpose of Section \ref{sec-cost functions} where these
costs functions are derived for a large family of random dynamics.
The proofs of our main results are done at Section
\ref{sec-proofs}. They are partly based on two
$\Gamma$-convergence results which are stated and proved at
Sections \ref{sec-Gamma1} and \ref{sec-Gamma2}. Finally, we recall
some basic notions about large deviation theory at the Appendix.

\subsection*{Notation} Let us introduce our main notations.

\textit{Measures}.\  For any topological space $X,$ we denote
$\mathrm{P}(X)$ the set of all Borel probability measures on $X$
and we endow it with the usual \emph{narrow  topology}
$\sigma(\mathrm{P}(X),\mathrm{C}_b(X))$ weakened by the space
$\mathrm{C}_b(X)$ of all continuous bounded functions on $X.$ We
also furnish $\mathrm{P}(X)$ with the corresponding Borel
$\sigma$-field.
\\
The push-forward of the measure $m$ by the measurable application
$f$ is denoted by $f\pf m$ and defined by $f\pf m(A)=m(f^{-1}(A))$
for any measurable set $A.$
\\
The Dirac measure at $a$ is denoted by $\delta_a.$

\textit{Measures on a path space}.\ We take a polish space $\XX$
which is furnished with its Borel $\sigma$-field. The relevant
space of  paths from the time interval $[0,1]$ to the state space
$\XX$ is either the space $\Omega=C([0,1],\XX)$ of all continuous
paths, or the space $\Omega=D([0,1],\XX)$ of paths which are left
continuous and right limited (c\`adl\`ag\footnote{This is the
french acronym for \textit{continu \`a droite et limit\'e \`a
gauche}.}) paths. We denote $X=(X_t)_{t\in[0,1]}$ the canonical
process which is defined for all $t\in[0,1]$ by
$$
X_t(\omega):=\omega_t,\quad
\omega=(\omega_t)_{t\in[0,1]}\in\Omega.
$$
For each $t\in\ii,$ $X_t$ is the position at time $t$ which is
seen as an application on $\Omega.$ Of course, $X$ is the identity
on $\Omega.$ The set $\Omega$ is endowed with the $\sigma$-field
$\sigma(X_t,t\in\ii)$ which is generated by the canonical process.
It is known that it matches with the Borel $\sigma$-field of
$\Omega$ when $\Omega$ is furnished with the Skorokhod
topology\footnote{In the special case where the paths are
continuous: $\Omega=C([0,1],\XX)$, this topology reduces to the
topology of uniform convergence.} which turns $\Omega$ into a
polish space, see \cite{Bil68}.  We denote $\PX,$ $\PXX$ and $\PO$
the set of all probability measures on $\XX,$ $\XXX=\XX\times\XX$
and $\Omega$ respectively. For any $P\in\PO,$ i.e.\ $P$ is the law
of a random path, we denote
    $$
P_t:=(X_t)\pf P\in \PX,\quad t\in\ii.
    $$
In particular, $P_0$ and $P_1$ are the laws of the initial and
final random positions under $P.$ Also useful is the joint law of
the initial and final positions
    $$
\bas P01:=(X_0,X_1)\pf P\in\PXX.
    $$
Of course, $\bas P01$ carries more information than the couple
$(P_0,P_1)$ because of the correlation structure. A similar remark
holds for $P\in\PO$ and $(P_t;t\in\ii)\in\PX^{\ii}.$ We denote the
disintegration of $P$ with respect to $(X_0,X_1):$
    $
P(d\omega)=\IXX P^{xy}(d\omega)\,\bas P01(dxdy)
    $
where
$$
    \haut Pxy(\cdot):=P(\cdot\mid X_0=x,X_1=y),\quad x,y\in\XX
$$
is the conditional law of $X$ knowing that $X_0=x$ and $X_1=y$
under $P.$ Its is usually called the bridge of $P$ between $x$ and
$y.$
\\
When working with the product space $\XXX,$ one sees the first and
second factors $\XX$ as the sets of initial and final states
respectively. Therefore, the canonical projections are denoted
$X_0(x,y):=x$ and $X_1(x,y):=y,$ $(x,y)\in\XXX.$ We denote the
marginals of the probability measure $\pi\in\PXX$ by
$\pi_0:=(X_0)\pf\pi\in\PX$ and $\pi_1:=(X_1)\pf\pi\in\PX.$

\textit{Functions}.\ Recall that a function
$f:X\to(-\infty,\infty]$ is said to be \lsc\ on the topological
space $X$ if all its sublevel sets $\{f\le a\},$ $a\in\mathbb{R}$
are closed. It is said to be \emph{coercive} if $X$ is assumed to
be Hausdorff and its sublevel sets are compact.
\\
Let $X$ and $Y$ be two topological vector spaces equipped with a
duality bracket $\langle x,y\rangle\in\mathbb{R},$ that is a
bilinear form on $X\times Y.$ The convex conjugate $f^*$ of
$f:X\to (-\infty,\infty]$ with respect to this duality bracket is
defined by
$$
    f^*(y):=\sup_{x\in X}\{\langle x,y\rangle
    -f(x)\}\in[-\infty,\infty],
    \quad y\in Y.
$$
It is a convex $\sigma(Y,X)$-\lsc\ function.
\\
The relative entropy of the probability $P$ with respect to the
probability $R$ is
$$
    H(P|R):=\left\{\begin{array}{ll}
      \int \log\left(\frac{dP}{dR}\right)\,dP\in[0,\infty] & \textrm{if }P\ll R \\
      \infty & \textrm{otherwise.} \\
    \end{array}\right.
$$

\section{Monge-Kantorovich and Schr\"odinger problems}\label{sec-problems}

In this section we present the Monge-Kantorovich optimal transport
problem  and the Schr\"odinger entropy minimization problem. Then,
we show informally with the aid of Schr\"odinger's original
example that they are connected to each other, by letting some
fluctuation coefficient tend to zero.
\\
\emph{Warning}. This informal section contains probabilistic
material. The probability-allergic reader can skip it without
harm. Nevertheless, it also contains some very clever ideas of
Schr\"odinger which acted as a guide for the author.

\subsection*{The Monge-Kantorovich optimal transport problem} Let
$c:\XXX\to[0,\infty]$ be a \lsc\ function on $\XXX$ with possibly
infinite values. For any $x,y\in\XX,$ $c(x,y)$ is interpreted as
the cost for transporting a unit mass from  location $x$ to
location $y.$ Let $\mu_0,\mu_1\in\PX$ be two prescribed
probability measures on $\XX.$ An admissible transport plan from
$\mu_0$ to $\mu_1$ is any probability measure $\pi\in\PXX$ which
has its first and second marginals equal to $\pi_0=\mu_0$ and
$\pi_1=\mu_1,$ respectively. For such a $\pi,$
$$
    \IXX c\,d\pi\in[0,\infty]
$$
is interpreted as the total cost for transporting $\mu_0$ to
$\mu_1$ when choosing the plan $\pi.$ The Monge-Kantorovich
optimal transport problem is the corresponding minimization
problem, i.e.\
\begin{equation}\label{MK}
     \IXX c\,d\pi \rightarrow \min;\quad \pi\in\PXX : \pi_0=\mu_0,
    \pi_1=\mu_1 \tag{MK}.
\end{equation}
A minimizer $\pih\in \PXX$ is called an optimal plan and the
minimal value $\inf$\eqref{MK} $\in[0,\infty]$ is the optimal
transport cost. Remark that \eqref{MK} is a convex minimization
problem. But, as it is \emph{not} a \emph{strictly} convex
problem, it might admit several solutions.

\subsection*{The Schr\"odinger entropy minimization problem}

Take a reference process $R$ on $\Omega$ (the unusual letter $R$
is chosen as a reminder of \emph{\textbf{r}eference}). By this, it
is meant a positive $\sigma$-finite measure on $\Omega$ which is
not necessarily bounded. Consider $n$ \emph{independent} random
dynamic particles $(Y^i; 1\le i\le n)$ where each random
realization of $Y^i$ lives in $\Omega.$ More specifically, $(Y^i;
1\le i\le n)$ is a collection of independent random paths where
for each $i,$ the law of $Y^i$ is
\begin{equation}\label{eq-12}
    \textrm{Law}(Y^i\mid Y^i_0)=R(\cdot\mid X_0=Y^i_0)\in\PO
\end{equation}
and  $(Y^i_0; 1\le i\le n)$ should be interpreted as the random
initial positions.

\begin{example}[Schr\"odinger's heat bath]\label{ex-01}
As a typical example, one can take $R$ to be the law of the
Brownian motion (Wiener process) on $\XX=\mathbb{R}^d$ with
diffusion coefficient $\sigma^2$ and the Lebesgue measure as its
initial distribution. The random motions are described by
$$
    Y^i_t=Y^i_0+\sigma B^i_t,\quad 1\le i \le n, t\in\ii
$$
where $Y^1_0,\dots,Y^n_0$ are random independent initial
positions, $B^1,\dots,B^n$ are independent Brownian motions with
initial position $B^i_0=0\in\mathbb{R}^d,$ $i=1,\dots,n,$ and
$\sigma>0$ is the square root of the temperature $\sigma^2$.
\end{example}

Schr\"odinger's original problem \cite{Sch32} which is based on
this specific example can be stated as follows. Suppose that at
time $t=0$ you observe a large particle system approximately in
the configuration $\mu_0\in \mathrm{P}(\mathbb{R}^d).$ The law of
large numbers tells you that with a very high probability you
observe the system at time $t=1$ in a configuration very near the
convolution $\mu_0
* \mathcal{N}(0,\sigma^2\mathrm{I})$ of $\mu_0$ and the normal
density with mean 0 and covariance matrix $\sigma^2\mathrm{Id}$ in
$\mathbb{R}^d.$ But, since the number $n$ of particle is
\emph{finite}, it is still possible (with a tiny probability of
order $e^{-an}$ with $a>0$) to observe the system at time $t=1$ in
a configuration which is significantly far form the expected
profile of distribution $\mu_0 *
\mathcal{N}(0,\sigma^2\mathrm{I}).$ Now, Schr\"odinger's question
is\footnote{Schr\"odinger's french words are: ``Imaginez que vous
observez un syst\`eme de particules en diffusion, qui soient en
\'equilibre thermodynamique. Admettons qu'\`a un instant donn\'e
$t_0$ vous les ayez trouv\'e en r\'epartition \`a peu pr\`es
uniforme et qu'\`a $t_1>t_0$ vous ayez trouv\'e un \'ecart
spontan\'e et \emph{consid\'erable} par rapport \`a cette
uniformit\'e. On vous demande de quelle mani\`ere cet \'ecart
s'est produit. Quelle en est la mani\`ere la plus probable ?"}:
\emph{Suppose that you observe the system at time $t=1$ in a
configuration which is approximately $\mu_1\in
\mathrm{P}(\mathbb{R}^d)$ and that $\mu_1$ is significantly
different from $\mu_0*\mathcal{N}(0,\sigma^2\mathrm{I}),$ what is
the most likely path of the whole system from $\mu_0$ to $\mu_1$
during the time interval $\ii?$}  In \cite{Sch32}, Schr\"odinger
gave the complete answer to this question with a proof based on
Stirling's formula. Although proved informally, there is nothing
significant to be added today to his answer.

The modern way of addressing this problem is in terms of large
deviations, see \cite{DZ} for an excellent overview of the large
deviation theory (a short reminder about large deviation theory is
also given at the Appendix). This has been done by F\"ollmer in
his Saint-Flour lecture notes \cite{Foe85}. Denoting $\delta_a$
the unit mass Dirac measure at $a,$ the whole system is described
by its empirical measure
$$
    L^n:=\frac 1n\sum_{i=1}^n\delta_{Y^i}\in \PO.
$$
It is a $\PO$-valued random variable\footnote{Strictly speaking
$L^n$ is a measurable function with its values in $\PO$ and the
statement ``$L^n\in\PO$'' is an abuse of notation. Nevertheless it
is a useful shorthand which will be used below without warning.}
which contains all the information about the dynamic system up to
any permutation of the labels of the particles. Nothing is lost
when the particles are indistinguishable. It also contains more
information than the random path
$$
t\in\ii\mapsto L^n_t:=\frac 1n\sum_{i=1}^n\delta_{Y^i_t}\in \PX
$$
which describes the evolution of the configurations. The observed
initial and final configurations are the empirical measures
    $$
    L^n_0=\frac 1n\sum_{i=1}^n\delta_{Y^i_0},\
    L^n_1=\frac 1n\sum_{i=1}^n\delta_{Y^i_1}\in\PX.
$$
Now, we give an informal presentation of the answer to
Schr\"odinger's question. For a rigorous treatment, one can have a
look at \cite{Foe85}. Take $\mu_0,\mu_1\in\PX$ and
$C_0^\delta,C_1^\delta$ two $\delta$-neighborhoods in $\PX$ (with
respect to a distance on $\PX$ compatible with the narrow topology
$\sigma(\PX,\CX)$) of $\mu_0$ and $\mu_1$ respectively. One can
recast Schr\"odinger's problem as follows: Find the measurable
sets $A\subset \PO$ such that the conditional probability
$$
\mathbb{P}(L^n\in A\mid L^n_0\in C^\delta_0,L^n_1\in C^\delta_1)
:=\frac{\mathbb{P}(L^n\in A,L^n_0\in C^\delta_0,L^n_1\in
C^\delta_1)}{\mathbb{P}(L^n_0\in C^\delta_0,L^n_1\in C^\delta_1)}
$$
is maximal when $n$ is large and $\delta$ is small. We introduced
the $\delta$-blowups $C_0^\delta$ and $C_1^\delta$ to prevent from
dividing by zero. A slight variant of Sanov's theorem\footnote{If
$R$ is a \emph{probability} measure, then this is really Sanov's
theorem and $H(P|R)\in[0,\infty].$}, a standard large deviation
result whose exact statement is in terms of large deviation
principle, see Theorem \ref{res-01} at the Appendix, states that
\begin{equation*}
    \mathbb{P}(L^n\in A)\underset{n\rightarrow\infty}\asymp
    \exp\left(-n \inf\left\{H(P|R);P\in A\right\}\right),\quad
    A\subset\PO
\end{equation*}
for any measurable subset $A$ of $\PO$ (by the way, one must
define a $\sigma$-field on $\PO$) where
$$
    H(P|R):=\IO
    \log\left(\frac{dP}{dR}\right)\,dP\in(-\infty,\infty],
    \quad P\in\PO
$$
is the relative entropy of $P$ with respect to $R$. Under
integrability conditions on $\mu_0$ and $\mu_1$ which insure that
there exists some $P\in\PO$ such that $P_0=\mu_0,$ $P_1=\mu_1$ and
$H(P|R)<\infty,$ one deduces that
\begin{equation*}
      \mathbb{P}(L^n\in A\mid L^n_0\in C^\delta_0,L^n_1\in C^\delta_1)
      \underset{n\rightarrow\infty}\asymp
      \frac{\exp\left(-n \inf\left\{H(P|R);P:P\in A,P_0\in C_0^\delta,P_1\in C_1^\delta \right\}\right)}
      {\exp\left(-n \inf\left\{H(P|R); P: P_0\in C_0^\delta,P_1\in
      C_1^\delta\right\}\right)}.
\end{equation*}
It follows that we have the conditional law of large numbers
$$
    \Lim n \mathbb{P}(L^n\in A\mid L^n_0\in C^\delta_0,L^n_1\in C^\delta_1)
    =\left\{\begin{array}{ll}
      1 & \textrm{if }A\ni \Ph^{\delta} \\
      0 & \textrm{otherwise} \\
    \end{array}\right.
$$
where $\Ph^{\delta}$ is the unique solution of the entropy
minimization problem
\begin{equation*}
    H(P|R) \rightarrow \min; \quad P\in\PO: P_0\in C_0^\delta, P_1\in C_1^\delta.
\end{equation*}
The uniqueness  comes directly from the \emph{strict} convexity of
the relative entropy $H(\cdot|R).$ We finally see that, under some
conditions on the limits $C_0^\delta
\underset{\delta\downarrow0}\rightarrow \{\mu_0\}$ and $C_1^\delta
\underset{\delta\downarrow0}\rightarrow \{\mu_1\},$ the solution
to the Schr\"odinger problem is
$$
   \lim_{\delta\downarrow0} \Lim n \mathbb{P}(L^n\in A\mid L^n_0\in C^\delta_0,L^n_1\in C^\delta_1)
    =\left\{\begin{array}{ll}
      1 & \textrm{if }A\ni \Ph \\
      0 & \textrm{otherwise} \\
    \end{array}\right.
$$
where $\Ph$ is the unique solution of the \emph{Schr\"odinger
entropy minimization problem}:
\begin{equation}\label{S}
    H(P|R)\rightarrow \min;\quad  P\in\PO: P_0=\mu_0, P_1=\mu_1 \tag{S} .
\end{equation}
If one prefers relative entropies with respect to probability
measures, $\Ph$ is also the unique solution to
\begin{equation}\label{eq-01}
    H(P|R^{\mu_0})\rightarrow\min;\quad P\in\PO: P_0=\mu_0, P_1=\mu_1
\end{equation}
where the $\sigma$-finite measure $R$ has been replaced by the
probability measure
\begin{equation}\label{eq-13}
     R^{\mu_0}(d\omega):=\IX R(d\omega\mid X_0=x)\,\mu_0(dx)\in\PO
\end{equation}
which is the law of the process with initial distribution
$\mu_0\in\PX$ and the same dynamics as $R.$ This last formulation
is analytically suitable, but it introduces an artificial time
asymmetry. Nevertheless, we keep it because it will be useful.
Remark that
$$
    H(P|R^{\mu_0})\in[0,\infty]
$$
is nonnegative and the minimization problems \eqref{S} and
\eqref{eq-01} share the same minimizer under the constraint
$P_0=\mu_0$ since
$H(P|R^{\mu_0})=H(P|R)-\IX\frac{d\mu_0}{dR_0}\,dP_0=H(P|R)-\IX\frac{d\mu_0}{dR_0}\,d\mu_0$.

The minimization problem \eqref{S} looks like the
Monge-Kantorovich problem \eqref{MK}, but we can do better in this
direction, relying on the tensorization property of the relative
entropy. Namely, for any measurable function
$\Phi:\Omega\to\mathcal{Z}$ where $\mathcal{Z}$ is any polish
space with its Borel $\sigma$-field, we have
\begin{equation}\label{eq-20}
    H(P|R)=H(\phi\pf P|\phi\pf R)+\int_{\mathcal{Z}}
H\Big(P(\cdot\mid\phi=z)\Big|R(\cdot\mid\phi=z)\Big)\,\phi\pf
P(dz).
\end{equation}
With $\phi=(X_0,X_1),$ this gives us
\begin{equation}\label{eq-03}
    H(P|R)=H(\bas P01|\bas R01)+\IXX H\Big(\haut Pxy \Big|\haut
Rxy\Big)\,\bas P01(dxdy).
\end{equation}
Now, decomposing the marginal constraint $P_0=\mu_0,$ $P_1=\mu_1$
into $\bas P01=\pi\in\PXX$ and $(X_0)\pf\pi=\mu_0,$
$(X_1)\pf\pi=\mu_1$ we obtain
\begin{eqnarray*}
    &&\inf\{H(P|R); P\in\PO:P_0=\mu_0,P_1=\mu_1\}\\
    &=&\inf\{\inf[H(P|R);P\in\PO:\bas P01=\pi];
    \pi\in\PXX:\pi_0=\mu_0,\pi_1=\mu_1\}.
\end{eqnarray*}
With \eqref{eq-03}, we see that the inner term is
\begin{eqnarray*}
    &&\inf[H(P|R);P\in\PO:\bas P01=\pi]\\
    &=& H(\pi|\bas R01)+\inf\left[\IXX H\Big(\haut Pxy \Big|\haut Rxy\Big)\,\pi(dxdy);P\in\PO:\bas
    P01=\pi\right]\\
    &=& H(\pi|\bas R01)
\end{eqnarray*}
where the $\inf$ is uniquely attained when $\haut Pxy=\haut Rxy,$
for $\pi$-almost every $(x,y),$ since in this case $0=H\Big(\haut
Pxy \Big|\haut Rxy\Big)$ which is the minimal value of the
relative entropy. This also shows that for each $\pi\in\PXX,$
\begin{equation}\label{eq-24}
    \inf[H(P|R);P\in\PO:\bas P01=\pi]=H(R^\pi|R)=H(\pi|\bas R01)
\end{equation}
where
\begin{equation}\label{eq-06}
    R^\pi(\cdot):=\IXX \haut Rxy(\cdot)\,\pi(dxdy)
\end{equation}
is the mixture of the bridges $\haut Rxy$ with $\pi$ as a mixing
measure. Hence, the solution of \eqref{S} is
$$
\Ph=R^{\pih}
$$
where $\pih\in\PXX$ is the unique solution of
\begin{equation}\label{eq-04}
    H(\pi|\bas R01)\rightarrow\min;\quad  \pi\in\PXX: \pi_0=\mu_0, \pi_1=\mu_1.
\end{equation}

\subsection*{Connecting \eqref{S} and \eqref{MK}}
The problem \eqref{eq-04} is similar to \eqref{MK}, but it remains
something to do in order to connect $\bas R01$ with some cost
function $c.$ We are going to do it in the special case which is
described at Example \ref{ex-01}, by letting the temperature
$$
1/k:=\sigma^2
$$
tend to zero where $k\ge1$ describes the positive integers. The
general situation will be investigated later at Section
\ref{sec-results}.
\\
Let us make the $k$-dependence explicit in our notation. We denote
$R^k$ the law of the process $Y^k$ which is defined by
\begin{equation}\label{eq-08}
    Y^k_t=Y_0+\sqrt{1/k} B_t,\quad 0\le t\le1
\end{equation}
with $Y_0$ having the Lebesgue measure as its distribution. In
particular, the joint law of the initial and final positions under
this reference process at positive temperature $1/k$ is
$$
    \bas {R^k}01(dxdy)
    =dx\,
    (2\pi/k)^{-d/2}\exp\left(-k\frac{|y-x|^2}{2}\right)\,dy.
$$
Rewriting \eqref{eq-04} with the $k$-dependence made explicit, we
get
\begin{equation}\label{eq-05}
    \frac 1k H(\pi|\bas {R^k}01)\rightarrow\min;\quad \pi\in\PXX:\pi_0=\mu_0, \pi_1=\mu_1.  \tag{${\widetilde{\mathcal{S}}}_{01}^k$}
\end{equation}
We ``tilde'' the name of this problem, because there will be
another ``untilded'' problem later:
\begin{equation}\label{S01e}
   \frac 1k H(\pi|\bas {R^{\mu_0,k}}01)\rightarrow\min;\quad
   \pi\in\PXX:\pi_0=\mu_0, \pi_1=\mu_1^k,  \tag{${\mathcal{S}}_{01}^k$}
\end{equation}
with better convergence properties, where the constraint $\mu_1$
is replaced by the ``moving'' constraint $\mu_1^k$ which is
indexed by $k$ and satisfies $\Lim k\mu_1^k=\mu_1$ and $\bas
{R^k}01$ is replaced by $\bas {R^{k,\mu_0}}01,$ see \eqref{eq-13}.
\\
We have also introduced the renormalization $\frac 1k H(\cdot|\bas
{R^k}01).$ To see that $\frac 1k$ is the right multiplying factor,
suppose that $\pi$ is such that $H(\pi|\bas {R^k}01)<\infty.$ This
implies that $\pi\ll\bas {R^k}01\ll \lambda$ where
$\lambda(dxdy)=dxdy$ stands for the Lebesgue measure  on
$\mathbb{R}^d\times\mathbb{R}^d.$ We see that
    $$
H(\pi|\bas {R^k}01)=
   H(\pi|\lambda) +\frac d2\log(2\pi/k)
    +k\int_{\mathbb{R}^d\times\mathbb{R}^d} \frac{|y-x|^2}2\,\pi(dxdy).
    $$
If we assume that $\pi$ satisfies
$\int_{\mathbb{R}^d\times\mathbb{R}^d}
\frac{|y-x|^2}2\,\pi(dxdy)<\infty,$ we obtain that the Boltzmann
entropy $H(\pi|\lambda)$ is finite and
$$
    \Lim k \frac1k H(\pi|\bas {R^k}01)=\int_{\mathbb{R}^d\times\mathbb{R}^d}\frac{|y-x|^2}2\,\pi(dxdy)
$$
which is the cost for transporting $\pi_0$ to $\pi_1$ with respect
to the quadratic cost
$$
    c(x,y)=|y-x|^2/2,\quad x,y\in\mathbb{R}^d.
$$
\textit{This indicates that \eqref{S01e} might converge to
\eqref{MK} in some sense, as $k$ tends to infinity.} Indeed, this
will be made precise and proved in the subsequent pages.
\\
The renormalized problem \eqref{S} with the dependence on $k$ made
explicit is
\begin{equation}\label{eq-07}
     \frac1k H(P|R^k)\rightarrow\min;\quad  P\in\PO:P_0=\mu_0, P_1=\mu_1. \tag{${\widetilde{\mathcal{S}}}^k$}
\end{equation}
Because of \eqref{eq-06}, its solution is $R^{k,\tilde{\pi}^{k}}$
where $\tilde{\pi}^{k}$ is the solution to \eqref{eq-05}. Remark
that for two distinct $k, k'>0,$ the supports of $R^k$ and
$R^{k'}$ are disjoint subsets of $\Omega.$ Therefore, at the level
of the process laws $P\in\PO,$ we see that for all $P\in\PO,$
$H(P|R^k)=\infty$ for every $k\ge1$ except possibly one. It
appears that the pointwise limit of $\frac1k H(\cdot|R^k)$ as $k$
tends to infinity is irrelevant. We shall see that the good notion
of convergence is that of $\Gamma$-convergence. Also, we shall
need the following ``untilded'' variant of \eqref{eq-07}:
\begin{equation}\label{Se}
     \frac1k H(P|R^{k,\mu_0})\rightarrow\min;\quad P\in\PO: P_0=\mu_0, P_1=\mu_1^k,
\tag{${\mathcal{S}}^k$}
\end{equation}
where $\Lim k\mu_1^k=\mu_1.$

\section{Statement of the main results}\label{sec-results}

The statements of our results is in terms of $\Gamma$-convergence
and large deviation principle. We start introducing their
definitions.

\subsection*{$\Gamma$-convergence}

We refer to the monograph by Dal~Maso \cite{DalMaso} for a clear
exposition of the subject. Recall that if it exists, the
$\Gamma$-limit of the sequence $\seq fk$ of
$(-\infty,\infty]$-valued functions on a topological space $X$ is
given for all $x$ in $X$ by
\begin{equation*}
    \Glim k f_k(x)=\sup_{V\in \mathcal{N}(x)}\Lim k \inf_{y\in V}
    f_k(y)
\end{equation*}
where $\mathcal{N}(x)$ is the set of all neighborhoods of $x.$ In
a metric space $X,$ this is equivalent to:
\begin{enumerate}[(i)]
    \item For any sequence $\seq xk$ such that $\Lim kx_k=x,$
\begin{equation*}
    \Liminf kf_k(x_k)\ge f(x)
\end{equation*}
 \item and  there exits a sequence $\seq{\tilde{x}}k$ such that $\Lim
 k\tilde{x}_k=x$ and
\begin{equation*}
    \Liminf kf_k(\tilde{x}_k)\le f(x).
\end{equation*}
\end{enumerate}
Item (i) is called the \emph{lower bound} and the sequence
$\seq{\tilde{x}}k$ in item (ii) is the \emph{recovery sequence}.

\subsection*{Large deviation principle}

We refer to the monograph by Dembo and Zeitouni \cite{DZ} for a
clear exposition of the subject. Let $X$ be a polish space
furnished with its Borel $\sigma$-field. One says that the
sequence $\seq{\gamma}n$ of probability measures on $X$ satisfies
the large deviation principle (LDP for short) with scale $n$ and
rate function $I,$ if for each Borel measurable subset $A$ of $X$
we have
\begin{equation}\label{eq-02}
    -\inf_{x\in\inter A}I(x)\overset{(\textrm{i})}\le\Liminf n \frac1n\log\gamma_n(A)
    \le \Limsup n \frac1n\log\gamma_n(A)\overset{(\textrm{ii})}\le-\inf_{x\in\cl A}I(x)
\end{equation}
where $\inter A$ and $\cl A$ are respectively the topological
interior and closure of $A$ in $X$ and the rate function
$I:X\to[0,\infty]$ is \lsc. The inequalities (i) and (ii) are
called respectively the \emph{LD lower bound} and \emph{LD upper
bound}, where LD is an abbreviation for large deviation. The LDP
is the exact statement of what was meant in previous section when
writing
$$
\gamma_n(A)\underset{n\rightarrow\infty}\asymp
\exp\left(-n\inf_{x\in A}I(x)\right)
$$
for ``all" $A\subset X.$

\subsection*{The main results}

For any topological space $X,$ we denote $\mathrm{P}(X)$ the set
of all Borel probability measures on $X$ and we endow it with the
usual weak topology $\sigma(\mathrm{P}(X),\mathrm{C}_b(X))$
weakened by the space $\mathrm{C}_b(X)$ of all continuous bounded
functions on $X.$ We also furnish $\mathrm{P}(X)$ with the
corresponding Borel $\sigma$-field.
\\
One says that a function $f:X\to(-\infty,\infty]$ is
\emph{coercive} if for any real $a\ge\inf f,$ the sublevel set
$\{f\le a\}$ is compact. This implies that $f$ is \lsc\ if $X$ is
Hausdorff, which will be the case of all the topological spaces in
the sequel.
\\
The convex analysis indicator of a set $A\subset X$ is defined by
$$
   \iota_{\{x\in A\}}= \iota_A(x)=\left\{\begin{array}{ll}
      0 & \textrm{if }x\in A \\
      \infty & \textrm{otherwise} \\
    \end{array}\right.,
    \quad x\in X.
$$
We keep the notation of Section \ref{sec-problems}. In particular,
$\XX$ is a polish space (metric complete and separable) with its
Borel $\sigma$-field and $$\Omega=D([0,1],\XX)$$ is the set of all
c\`adl\`ag $\XX$-valued paths endowed with the Skorokhod metric,
which turns it into a polish space.

\subsubsection*{Static version}
For each integer $k\ge1,$ we take a measurable kernel
$$(\pont \rho kx\in \PX;x\in\XX)$$ of probability
measures on $\XX$. We also take $\mu_0\in\PX,$  denote
$$
    \pont \rho{k}{\mu_0}(dxdy):=\mu_0(dx)\pont \rho kx(dy)\in\PXX
$$
and define the functions
$$
    \bas{\mathcal{C}}01^{k,\mu_0}(\pi):=\frac1k H(\pi| \pont \rho{k}{\mu_0})+\iota_{\{\pi_0=\mu_0\}},\ k\ge1;\quad
    \bas{\mathcal{C}}01^{\mu_0}(\pi):=\IXX c\,d\pi+\iota_{\{\pi_0=\mu_0\}},\quad
    \pi\in\PXX.
$$

\begin{proposition}\label{res-05}
We assume that for each $x\in\XX,$ the sequence $((X_1)\pf\pont
\rho kx)_{k\ge1}$ satisfies the LDP in $\XX$ with scale $k$ and
the \emph{coercive} rate function
$$
    c(x,\cdot):\XX\to[0,\infty]
$$
where $c:\XXX\to[0,\infty]$ is a \lsc\ function.
\\
Then, for any $\mu_0\in\PX$ we have:
    $
    \Glim k   \bas{\mathcal{C}}01^{k,\mu_0}=\bas{\mathcal{C}}01^{\mu_0}
    $
in $\PXX.$
\end{proposition}
Let us define the functions
\begin{eqnarray*}
  \bas T01^k(\mu_0,\nu)
  &:=&\inf\left\{\frac1k H(\pi| \pont
  \rho{k}{\mu_0});\pi\in\PXX:\pi_0=\mu_0,\pi_1=\nu\right\}\\
  &=&\inf\{\bas{\mathcal{C}}01^k(P);\pi\in\PXX:\pi_1=\nu\}
  ,\qquad \nu\in\PX
\end{eqnarray*}
and
\begin{eqnarray*}
   \bas T01(\mu_0,\nu)
   &:=&\inf\left\{\IXX c\,d\pi;\pi\in\PXX:\pi_0=\mu_0,\pi_1=\nu\right\}\\
   &=&\inf\{\bas{\mathcal{C}}01^{\mu_0}(\pi);\pi\in\PXX:\pi_1=\nu\},\qquad \nu\in\PX.
\end{eqnarray*}
 The subsequent result follows easily from Proposition
\ref{res-05}.
\begin{corollary}\label{res-06}
Under the assumptions of Proposition \ref{res-05}, for any
$\mu_0\in\PX$ we have
$$
    \Glim k  \bas T01^k(\mu_0,\cdot)=\bas T01(\mu_0,\cdot)
$$
on $\PX.$ In particular, for any $\mu_1\in\PX,$ there exists a
sequence $\Seq{\mu_1}k$ such that $\Lim k \mu_1^k=\mu_1$ in $\PX$
and $\Lim k \bas T01^k(\mu_0,\mu_1^k)=\bas
T01(\mu_0,\mu_1)\in[0,\infty].$
\end{corollary}

Now, let us consider a sequence of minimization problems  which is
a generalization of $\eqref{S01e}_{k\ge1}$ at Section
\ref{sec-problems}. It is given for each $k\ge1,$ by
\begin{equation}\label{S01k}
     \frac1k H(\pi|\pont \rho{k}{\mu_0})\rightarrow\min;
     \quad \pi\in\PXX:\pi_0=\mu_0,\pi_1=\mu_1^k
     \tag{S$_{01}^k$}
\end{equation}
where $\Seq{\mu_1}k$ is a sequence in $\PX$ as in Corollary
\ref{res-06}. The Monge-Kantorovich problem associated with
$\eqref{S01k}_{k\ge1}$ is
\begin{equation*}
    \IXX c\,d\pi\rightarrow\min;\quad \pi\in\PXX:\pi_0=\mu_0,\pi_1=\mu_1.
     \tag{MK}
\end{equation*}

The main result of the paper is the following theorem.

\begin{theorem}\label{res-07}
Under the assumptions of Proposition \ref{res-05}, for any
$\mu_0,\mu_1\in\PX$ we have
    $
    \Lim k\inf\eqref{S01k}=\inf\eqref{MK}\in[0,\infty].
    $
\\
Suppose that in addition $\inf\eqref{MK}<\infty,$ then for each
large enough $k,$ \eqref{S01k} admits a \emph{unique} solution
$\pih^k\in\PXX.$
\\
Moreover, any limit point of the sequence $\Seq{\pih}k$ in $\PXX$
is a solution to \eqref{MK}. In particular, if \eqref{MK} admits a
unique solution $\pih\in\PXX,$ then $\Lim k \pih^k=\pih$ in
$\PXX.$
\end{theorem}
Remark that $\Lim k\inf\eqref{S01k}=\inf\eqref{MK}$ is a
restatement of $\Lim k \bas T01^k(\mu_0,\mu_1^k)=\bas
T01(\mu_0,\mu_1)$ in Corollary \ref{res-06}.

Proposition \ref{res-05} and Theorem \ref{res-07} admit a dynamic
version.

\subsubsection*{Dynamical version}
For each integer $k\ge1,$ we take a measurable kernel
$$
    (\pont Rkx\in \mathrm{P}(\Omega^x);x\in\XX)
$$
of probability measures on $\Omega,$ with
$$
    \Omega^x:=\{X_0=x\}.
$$
We have in mind the situation where $R^k\in\PO$ is the law of
stochastic process and  $\pont Rkx=R^k(\cdot\mid X_0=x)$ is its
conditional law knowing that $X_0=x,$ see \eqref{eq-12}. For any
$\mu_0\in\PX,$ denote
$$
    \pont R{k}{\mu_0}(\cdot):=\IX \pont Rkx(\cdot)\,\mu_0(dx)\in\PO,
    \quad
    \pont {\bas R01}{k}{\mu_0}(\cdot):=\IX \pont {\bas
    R01}kx(\cdot)\,\mu_0(dx)\in\PXX.
$$
We see that  $\pont R{k}{\mu_0}$ is the law of a stochastic
process with initial law $\mu_0$ and its dynamics determined by
$(\pont Rkx;x\in\XX)$ where $x$ must be interpreted as an initial
position, while $\pont {\bas R01}{k}{\mu_0}=(X_0,X_1)\pf \pont
R{k}{\mu_0}$ is the joint law of the initial and final positions
under $\pont R{k}{\mu_0}$. Let us define the functions
$$
    \mathcal{C}^{k,\mu_0}(P):=\frac1k H(P|\pont R{k}{\mu_0})+\iota_{\{P_0=\mu_0\}},\ k\ge1,\quad
    \mathcal{C}^{\mu_0}(P):=\IO C\,dP+\iota_{\{P_0=\mu_0\}},\quad
    P\in\PO
$$
where $C:\Omega\to[0,\infty]$ is a \lsc\ function.

\begin{proposition}\label{res-02}
We assume that for each $x\in\XX,$ the sequence $(\pont
Rkx)_{k\ge1}$ satisfies the LDP in $\Omega$ with scale $k$ and the
\emph{coercive} rate function
$$
    C^x=C+\iota_{\{X_0=x\}}:\Omega\to[0,\infty]
$$
where $C:\Omega\to[0,\infty]$ is a \lsc\ function.
\\
Then, for any $\mu_0\in\PX$ we have:
    $
    \Glim k  \mathcal{C}^{k,\mu_0}=\mathcal{C}^{\mu_0}
    $
in $\PO.$
\end{proposition}

Let us define the functions
\begin{eqnarray*}
  T^k(\mu_0,\nu)
  &:=& \inf\left\{\frac1k H(P|\pont R{k}{\mu_0});P\in\PO:P_0=\mu_0,P_1=\nu\right\}\\
  &=&\inf\{\mathcal{C}^{k,\mu_0}(P);P\in\PO:P_1=\nu\},\quad\nu\in\PX
\end{eqnarray*}
and
\begin{eqnarray*}
  T(\mu_0,\nu)
  &:=& \inf\left\{\IO C\,dP;P\in\PO:P_0=\mu_0,P_1=\nu\right\}\\
  &=&\inf\{\mathcal{C}^{\mu_0}(P);P\in\PO:P_1=\nu\},\quad \nu\in\PX.
\end{eqnarray*}

\begin{corollary}\label{res-03}
Under the assumptions of Proposition \ref{res-02}, we have
$$
    \Glim k  T^k(\mu_0,\cdot)=T(\mu_0,\cdot)
$$
on $\PX.$ In particular, for any $\mu_1\in\PX,$ there exists a
sequence $\Seq{\mu_1}k$ such that $$\Lim k \mu_1^k=\mu_1$$ in
$\PX$ and $\Lim k T^k(\mu_0,\mu_1^k)=T(\mu_0,\mu_1)\in[0,\infty].$
\end{corollary}

Now, let us consider a sequence of minimization problems  which is
a generalization of $\eqref{Se}_{k\ge1}$ at Section
\ref{sec-problems}. It is given for each $k\ge1,$ by
\begin{equation}\label{Sk}
     \frac1k H(P|\pont R{k}{\mu_0})\rightarrow\min;
     \quad P\in\PO:P_0=\mu_0,P_1=\mu_1^k
     \tag{S$^k$}
\end{equation}
where $\Seq{\mu_1}k$ is a sequence in $\PX$ as in Corollary
\ref{res-03}. The dynamic Monge-Kantorovich problem associated
with  $\eqref{Sk}_{k\ge1}$ is
\begin{equation}\label{MKdyn}
    \IO C\,dP\rightarrow\min;\quad P\in\PO:P_0=\mu_0,P_1=\mu_1.
     \tag{MK$_{\textrm{dyn}}$}
\end{equation}
Indeed, next result states that $\eqref{Sk}_{k\ge1}$ tends to
\eqref{MKdyn} in the sense that not only the values of
$\eqref{Sk}_{k\ge1}$ tend to $\inf\eqref{MKdyn},$ but also the
minimizers of $\eqref{Sk}_{k\ge1}$ tend to some minimizers of the
limiting problem \eqref{MKdyn}.

\begin{theorem}\label{res-04}
Under the assumptions of Proposition \ref{res-02}, for any
$\mu_0,\mu_1\in\PX$ we have
    $
    \Lim k\inf\eqref{Sk}=\inf\eqref{MKdyn}\in[0,\infty].
    $
\\
Suppose that in addition $\inf\eqref{MKdyn}<\infty,$ then for each
large enough $k,$ \eqref{Sk} admits a \emph{unique} solution
$\Ph^k\in\PO.$
\\
Moreover, any limit point of the sequence $\Seq{\Ph}k$ in $\PO$ is
a solution to $\eqref{MKdyn}.$ In particular, if \eqref{MKdyn}
admits a unique solution $\Ph\in\PO,$ then $\Lim k \Ph^k=\Ph$ in
$\PO.$
\end{theorem}

\subsubsection*{From the dynamic to a static version}

Once we have the dynamic results, the static ones can be derived
by means of the continuous mapping $P\in\PO\mapsto (X_0,X_1)\pf
P=\bas P01\in\PXX.$ The LD tool which is behind this transfer is
the contraction principle which is recalled at Theorem
\ref{res-08} below. The connection between the dynamic cost $C$
and the static cost $c$ is
\begin{equation}\label{eq-10}
     c(x,y):=\inf\{C(\omega);\omega\in\Omega:\omega_0=x,\omega_1=y\}\in [0,\infty],
    \quad x,y\in\XXX.
\end{equation}
This identity is connected to the \emph{geodesic problem}:
\begin{equation}\label{eq-11}
    C(\omega)\rightarrow\min;\quad
    \omega\in\Omega: \omega_0=x,\omega_1=y.
    \tag{G$^{xy}$}
\end{equation}
Since $C^x$ is coercive for all $x\in\XX,$ there exists at least
one solution to this problem, called a \emph{geodesic path},
provided that its value $c(x,y)$ is finite.

The above static results hold true for any $[0,\infty]$-valued
function $c$ satisfying the assumptions of Proposition
\ref{res-05} even if it is not derived from a dynamic rate
function $C$ via the identity \eqref{eq-10}. Note also that the
coerciveness of $C^x$ for all $x\in\XX,$ implies that
$y\in\XX\mapsto c(x,y)$ is coercive (the sublevel sets of
$c(x,\cdot)$ are continuous projections of the sublevel sets of
$C^x$ which are assumed to be compact). Nevertheless, it is not
clear at first sight that $c$ is  jointly (on $\XXX$) measurable.
Next result tells us that it is jointly \lsc.
\\
The coerciveness of $C^x$ also guarantees that the set of all
geodesic paths from $x$ to $y:$
$$
    \haut\Gamma
xy:=\{\omega\in\OO; \omega_0=x,\omega_1=y,C(\omega)=c(x,y)\}
$$
is a compact subset of $\OO$ which is nonempty as soon as
$c(x,y)<\infty.$ In particular, it is a Borel measurable subset.

\begin{theorem}\label{res-09}
Suppose that the assumptions of Proposition \ref{res-02} are
satisfied.
\begin{enumerate}
    \item Then, not only the dynamic results Corollary \ref{res-03} and
Theorem \ref{res-04} are satisfied with the cost function $C,$ but
also the static results Proposition \ref{res-05}, Corollary
\ref{res-06} and Theorem \ref{res-07} hold with the cost function
$c$ which is derived from $C$ by means of \eqref{eq-10}.
It is also true that $c$ is \lsc\ and
$\inf\eqref{MKdyn}=\inf\eqref{MK}\in[0,\infty].$
\end{enumerate}
Suppose in addition that $\mu_0,\mu_1\in\PX$ satisfy
$\inf\eqref{MK}:=\bas T01(\mu_0,\mu_1)<\infty,$ so that both
\eqref{MK} and \eqref{MKdyn} admit a solution.
\begin{enumerate}

    \item[(2)]  Then, for all large enough $k\ge1,$
\eqref{S01k} and \eqref{Sk} admit respectively a unique solution
$\pih^k\in\PXX$ and $\Ph^k\in\PO.$ Furthermore,
$$
    \Ph^k=R^{k,\pih^k}:=\IXX R^{k,xy}(\cdot)\,\pih^k(dxdy)
$$
which means that $\Ph^k$ is the $\pih^k$-mixture of the bridges
$R^{k,xy}$ of $R^k.$
    \item[(3)] The sets of solutions to \eqref{MK} and
    \eqref{MKdyn} are nonempty convex compact subsets of $\PXX$
    and $\PO$ respectively.
    \\
A probability $\Ph\in\PO$ is a solution to \eqref{MKdyn} if and
only if $\bas\Ph 01$ is a solution to \eqref{MK} and
\begin{equation}\label{eq-34}
    \haut\Ph xy(\haut\Gamma xy)=1,
    \quad \forall (x,y)\in\XX,\ \bas\Ph01\textrm{-a.e.}
\end{equation}
In particular, if \eqref{MK} admits a \emph{unique} solution
$\pih\in\PXX$ and  for $\pih$-almost every $(x,y)\in\XXX,$ the
geodesic problem \eqref{eq-11} admits a \emph{unique} solution
$\haut\gamma xy\in\Omega.$ Then, \eqref{MKdyn} admits the unique
solution
\begin{equation*}
    \Ph=\IXX \delta_{\haut\gamma xy}\,\pih(dxdy)\in\PO
\end{equation*}
which is the $\pih$-mixture of the Dirac measures at the geodesics
$\haut \gamma xy$ and $$\Lim k\Ph^k=\Ph$$ in $\PO.$
\end{enumerate}
\end{theorem}

\begin{remarks}\
\begin{enumerate}
     \item Formula \eqref{eq-34} simply means that $\Ph$ only
     charges geodesic paths. But we didn't write $\Ph(\Gamma)=1$
     since it is not clear that the set $\Gamma:=\bigcup_{x,y\in\XX}\haut\Gamma
     xy$ of all geodesic paths is measurable.
    \item In case of uniqueness as in the last statement of this theorem, the marginal flow of $\Ph$ is
$$
     \mu_t:=\Ph_t=\IXX \delta_{\haut{\gamma_t}
    xy}\,\pih(dxdy)\in\PX,\quad t\in\ii.
$$
It is the \emph{displacement interpolation} between $\mu_0$ and
$\mu_1.$
    \\ As a consequence of the abstract
disintegration result of the probability measures on a polish
space, the kernel $(x,y)\mapsto\delta_{\haut\gamma xy}$ is
measurable. This also means that $(x,y)\mapsto\haut\gamma xy$ is
measurable.
    \item If no uniqueness requirement is verified, then
$$
     \mu_t=(X_t)\pf\Ph\in\PX,\quad t\in\ii
$$
is also a good candidate for being called a displacement
interpolation between $\mu_0$ and $\mu_1.$
    \item The problem of knowing if $\Seq\pih k$ converges even if
    \eqref{MK} admits several solutions is left open in this
    article. It might be possible that this holds true and that
    the entropy minimization approximation selects a ``viscosity
    solution'' of \eqref{MK}.
\end{enumerate}
\end{remarks}

\subsection*{Back to Schr\"odinger's heat bath}
We illustrate these general results by means of Example
\ref{ex-01}. A well-known LD result is about the large deviations
of the $\RR^d$-valued process which we have already met at
\eqref{eq-08} and is defined by
\begin{equation}\label{eq-25}
    Y^{k,x}_t=x+\sqrt{1/k} B_t,\quad 0\le t\le1,
\end{equation}
where the initial condition $Y^{k,x}_0=x$ is deterministic,
$B=(B_t)_{0\le t\le1}$ is the Wiener process on $\RR^d$ and we
decided to take $
    \sigma^2=1/k.
$ with $k\ge1$ an integer.

\begin{theorem}[Schilder's theorem]
The sequence of random processes $(Y^{k,x})_{k\ge1}$ satisfies the
LDP in $\Omega=C([0,1],\RR^d)$ equipped with the topology of
uniform convergence with scale $k$ and rate function
$$
    C^x(\omega)= \int_\ii \frac{|\dot\omega_t|^2}2\,dt\in[0,\infty],
    \quad \omega\in\Omega
$$
if $\omega_0=x$ and $\omega$ is an absolutely continuous path (its
derivative is denoted by $\dot\omega$) and $C^x(\omega)=\infty,$
otherwise.
\end{theorem}
\noindent For a proof, see \cite[Thm 5.2.3]{DZ}.

With our notation, this corresponds to
\begin{equation*}
     C(\omega)= \left\{
    \begin{array}{ll}
      \int_\ii \frac{|\dot\omega_t|^2}2\,dt\in[0,\infty] & \textrm{if }\omega\in\Oac  \\
      \infty & \textrm{otherwise} \\
    \end{array}
     \right.,
    \quad \omega\in\Omega
\end{equation*}
where $\Oac$ is the space of all absolutely continuous paths
$\omega:[0,1]\to\RR^d.$ By  Jensen's inequality,  \eqref{eq-10}
leads us to
$$
    c(x,y)=|y-x|^2/2,\quad x,y\in\RR^d
$$
which is the well-known quadratic transport cost. Let $\pont
Rkx\in\PO$ denote the law of $\pont Ykx.$ Then $\pont
R{k}{\mu_0}(\cdot)=\int_{\RR^d}\pont Rkx(\cdot)\,\mu_0(dx) \in\PO$
is the law of
$$
    Y^{k}_t=Y_0+\sqrt{1/k} B_t,\quad 0\le t\le1,
$$
with initial law: $\textrm{Law}(Y_0)=\mu_0\in\mathrm{P}(\RR^d).$
Also denote $\pont \rho{k}{\mu_0}=(X_0,X_1)\pf \pont
R{k}{\mu_0}\in \mathrm{P}(\RR^d\times\RR^d),$ i.e.\
$$
    \pont \rho{k}{\mu_0}(dxdy)=\mu_0(dx)(2\pi/k)^{-d/2}\exp\left(-k\frac{|y-x|^2}{2}\right)\,dy.
$$
The above results tell us that if there exists some
$\pi^*\in\mathrm{P}(\RR^d\times\RR^d)$ such that $\pi^*_0=\mu_0,$
$\pi^*_1=\mu_1$ and $\int_{\RR^d\times\RR^d}
|y-x|^2\,\pi^*(dxdy)<\infty,$ then $\bas T01(\mu_0,\mu_1)<\infty,$
there exists a sequence $\Seq{\mu_1}k$ such that $\Lim
k\mu_1^k=\mu_1$ in $\mathrm{P}(\RR^d)$ and for any large enough
$k\ge1,$ the entropy minimization problem
$$
    \frac1k H(\pi|\pont\rho{k}{\mu_0})\rightarrow\min;\quad
    \pi\in \mathrm{P}(\RR^d\times\RR^d),\pi_0=\mu_0,\pi_1=\mu_1^k
$$
admits a unique solution $\pih^k\in \mathrm{P}(\RR^d\times\RR^d),$
the sequence $\Seq{\pih}k$ admits at least a limit point in
$\mathrm{P}(\RR^d\times\RR^d)$ and any such limit point is a
solution of the Monge-Kantorovich quadratic transport problem
$$
    \int_{\RR^d\times\RR^d} \frac{|y-x|^2}2\,\pi(dxdy)\rightarrow\min;\quad
    \pi\in\mathrm{P}(\RR^d\times\RR^d),\pi_0=\mu_0,\pi_1=\mu_1.
$$
In addition, we have $\Lim k  \frac1k
H(\pih^k|\pont\rho{k}{\mu_0})=\bas T01(\mu_0,\mu_1).$
\\
Moreover, for all large enough $k\ge1,$ the corresponding dynamic
problem
$$
    \frac1k H(P|\pont R{k}{\mu_0})\rightarrow\min;\quad
    P\in\PO,P_0=\mu_0,P_1=\mu_1^k
$$
has a unique solution $\Ph^k\in\PO$ which is given by
    $$
    \Ph^k=R^{k,\pih^k}=\int_{\RR^d\times\RR^d}\pont Rk{xy}(\cdot)\,\pih^k(dxdy)\in\PO,
    $$
the sequence $\Seq{\Ph}k$ admits at least a limit point in $\PO$
and any such limit point is a solution of the dynamic
Monge-Kantorovich quadratic transport problem
$$
    \int_{\Oac} \left[ \int_\ii \frac{|\dot\omega_t|^2}2\,dt\right]\,P(d\omega)\rightarrow\min;\quad
    P\in \mathrm{P}(\Oac),P_0=\mu_0,P_1=\mu_1.
$$
In the case where $\mu_0$ or $\mu_1$ is absolutely continuous with
respect to the Lebesgue measure on $\RR^d,$ it is well known
\cite{Bre91,McC95} that \eqref{MK} admits a unique solution
$\pih.$ By Theorem \ref{res-09}, we obtain $\Lim k \Ph^k=\Ph$ in
$\PO$ where
$$
    \Ph(\cdot)=\int_{\RR^d\times\RR^d}\delta_{[t\mapsto(1-t)x+ty]}(\cdot)\,\pih(dxdy)\in\PO
$$
and the corresponding displacement interpolation is the marginal
flow of $\Ph$ which is
$$
      \mu_t(\cdot)=\Ph_t(\cdot)=\int_{\RR^d\times\RR^d}\delta_{(1-t)x+ty}(\cdot)\,\pih(dxdy)\in\mathrm{P}(\RR^d),
      \quad t\in\ii.
$$
Moreover, the marginal flow of $\Ph^k$ is
$$
      \mu_t^k(\cdot)=\Ph^k_t(\cdot)=\int_{\RR^d\times\RR^d}\pont{R_t}k{xy}(\cdot)\,\pih^k(dxdy)\in\mathrm{P}(\RR^d),
      \quad t\in\ii
$$
and for each $t\in\ii,$ $\Lim k\mu^k_t=\mu_t$ in
$\mathrm{P}(\RR^d).$

Mikami's paper \cite{Mikami04} is in the context of
Schr\"odinger's heat bath based on the Wiener process as above.
Although no relative entropy nor $\Gamma$-convergence enter the
statements of \cite{Mikami04}'s results,  some of the previous
results about Schr\"odinger's heat bath are close to the main
results of \cite{Mikami04} which are proved by means of cyclical
monotonicity with a stochastic optimal control point of view.
Theorems \ref{res-07}, \ref{res-04}, \ref{res-09} apply to a large
class of optimal transport costs, see Section \ref{sec-cost
functions}. They shed a new light on \cite{Mikami04}'s results and
extend them in the sense that the reference process $R$ is not
restricted to be the Wiener process and the LD principle which is
satisfied by $\Seq Rk$ is not restricted to the setting of
Schilder's theorem.

\section{From stochastic processes to transport cost functions}\label{sec-cost functions}

We have just seen that Schilder's theorem leads to the quadratic
cost function. The aim of this section is to present a series of
examples of LD sequences $\Seq Rk$ in $\PO$ which give rise to
cost functions $c$ on $\XXX.$

\subsection*{Simple random walks on $\RR^d$}
Instead of \eqref{eq-25}, let us consider
\begin{equation}\label{eq-26}
     \pont Ykx_t=x+W^k_t,\quad 0\le t\le1,
\end{equation}
where for each $k\ge1,$ $W^k$ is a random walk. The law of $\pont
Ykx$ is our $\pont Rkx\in\PO.$

To build these random walks, one needs a sequence of independent
copies $\seq Zm$ of a random variable $Z$ in $\RR^d.$ For each
integer $k\geq 1,$  $W^k$ is the rescaled random walk defined for
all $0\leq t\leq 1,$ by
\begin{equation}\label{eq-43}
 W^k_t=\frac{1}{k}\sum_{j=1}^{\lfloor kt\rfloor} Z_j
% +\left(t-\frac{\lfloor kt\rfloor}{k}\right)Z_{\lfloor kt\rfloor +1}
\end{equation}
where $\lfloor kt\rfloor$ is the integer part of $kt.$ This
sequence satisfies a LDP which is given by Mogulskii's theorem. As
a pretext to set some notations, we recall its statement. The
logarithm of the Laplace transform of the law
$m_Z\in\mathrm{P}(\RR^d)$ of $Z$ is $\log
\int_{\RR^d}e^{\zeta\cdot z}\,m_Z(dz).$ Its convex conjugate is
\begin{equation}\label{cc}
    c_Z(v):=\sup_{\zeta\in\RR^d}\left\{\zeta\cdot v-\log \int_{\RR^d}e^{\zeta\cdot z}\,m_Z(dz)\right\},\quad v\in\RR^d.
\end{equation}
One can prove, see \cite{DZ}, that $c_Z$ is a convex
$[0,\infty]$-valued function which attains its minimum value $0$
at $v=\mathbb{E}Z=\int_{\RR^d}z\,m_Z(dz).$ Moreover, the closure
of its effective domain $\cl \{c_Z<\infty\}$ is the closed convex
hull of the topological support $\supp m_Z$ of the probability
measure $m_Z.$ Under the assumption \eqref{aa} below, it is also
\emph{strictly} convex.
\\
For each initial value $x\in\RR^d,$ we define the action
functional
\begin{equation*}
    C_Z^x(\omega):=\left\{\begin{array}{ll}
      \Iii c_Z(\dot{\omega}_t)\,dt & \textrm{if $\omega\in \Oac$ and }\omega_0=x \\
      +\infty & \textrm{otherwise} \\
    \end{array}\right.,\quad \omega\in \OO.
\end{equation*}

\begin{theorem}[Mogulskii's theorem]
Under the assumption
\begin{equation}\label{aa}
 \int_{\RR^d}e^{\zeta\cdot z}\,m_Z(dz)<+\infty, \quad
\forall \zeta\in\RR^d,
\end{equation}
for each $x\in\RR^d$ the sequence $(\pont Rkx)_{k\ge1}$ of the
laws of $(\pont Ykx)_{k\ge1}$ specified by \eqref{eq-26} satisfies
the LDP in $\Omega=D(\ii,\RR^d),$ equipped with its natural
$\sigma$-field and the topology of uniform convergence, with scale
$k$ and the coercive rate function $C_Z^x.$
\end{theorem}
For a proof see \cite[Thm 5.1.2]{DZ}. This result corresponds to
our general setting with
\begin{equation}\label{eq-30}
    C(\omega)=C_Z(\omega):=\left\{\begin{array}{ll}
      \Iii c_Z(\dot{\omega}_t)\,dt & \textrm{if }\omega\in \Oac \\
      +\infty & \textrm{otherwise} \\
    \end{array}\right.,\quad \omega\in \OO.
\end{equation}
Since $c_Z$ is a strictly convex function, the geodesic problem
\eqref{eq-11} admits as unique solution the \emph{constant
velocity geodesic}
\begin{equation}\label{eq-31}
    \haut\sigma xy:t\in\ii\mapsto (1-t)x+ty\in\RR^d.
\end{equation}
Now, let us only consider the final position $$\pont
Ykx_1=x+\frac{1}{k}\sum_{j=1}^{k} Z_j.
$$
Denote $\pont\rho kx=(X_1)\pf\pont Rkx\in\PX$ the law of $\pont
Ykx_1$.  By the contraction principle, see Theorem \ref{res-08} at
the Appendix, one deduces immediately from Mogulskii's theorem the
simplest result of LD theory which is the Cram\'er theorem.

\begin{corollary}[A complicated version of Cram\'er's theorem]
Under the assumption \eqref{aa}, for each $x\in\RR^d$ the sequence
$(\pont \rho kx)_{k\ge1}$ of the laws of  $(\pont Ykx_1)_{k\ge1}$
satisfies the LDP in $\RR^d$ with scale $k$ and the coercive rate
function
$$
   y\in\XX\mapsto c_Z(y-x)\in[0,\infty]\quad y\in\XX
$$
where $c_Z$ is given at \eqref{cc}.
\\
Furthermore,
$c_Z(v)=\inf\{C_Z(\omega);\omega\in\OO:\omega_0=x,\omega_1=x+v\}$
for all $x,v\in\RR^d.$
\end{corollary}
Last identity is a simple consequence of Jensen's inequality which
also lead us to \eqref{eq-31} a few lines earlier. Cram\'er's
theorem corresponds to the case when $x=0$ and only the deviations
of $\pont Yk0_1=\frac{1}{k}\sum_{j=1}^{k} Z_j$ in $\RR^d$ are
considered.

\begin{theorem}[Cram\'er's theorem]
Under the assumption \eqref{aa}, the sequence
$(\frac{1}{k}\sum_{j=1}^{k} Z_j)_{k\ge1}$  satisfies the LDP in
$\RR^d$ with scale $k$ and the coercive rate function $c_Z$ given
at \eqref{cc}.
\end{theorem}
For a proof, see \cite[Thm 2.2.30]{DZ}.

We have just described a general procedure which converts the law
$m_Z\in\mathrm{P(\RR^d)}$ into the cost functions $C_Z$ and $c_Z.$
Here are some examples with explicit computations.

\begin{examples}\label{ex-03}
We recall some well-known examples of Cram\'er transform $c_Z.$
\begin{enumerate}
    \item To obtain the quadratic cost function $c_Z(v)=|v|^2/2,$ choose $Z$
    as a standard normal random vector in $\RR^d:$ $m_Z(dz)=(2\pi)^{-d/2}\exp(-|z|^2/2)\,dz.$
    \item Taking $Z$ such that
    $\mathbb{P}(Z=+1)=\mathbb{P}(Z=-1)=1/2,$ i.e.\ $m_Z=(\delta_{-1}+\delta_{+1})/2$ leads to\\
    $c_Z(v)=\left\{%
\begin{array}{ll}
    [(1+v)\log(1+v)+(1-v)\log(1-v)]/2, & \hbox{if }-1<v<+1 \\
    \log 2, & \hbox{if } v\in\{-1,+1\}\\
    +\infty, & \hbox{if } v\not\in [-1,+1]. \\
\end{array}%
\right.    $
    \item If $Z$ has an exponential law with expectation 1, i.e.\ $m_Z(dz)=\1_{\{z\ge0\}}e^{-z}\,dz,$ then $c_Z(v)=v-1-\log
    v$ if $v>0$ and $c_Z(v)=+\infty$ if $v\leq 0.$
    \item If $Z$ has a Poisson law with expectation 1, i.e.\ $m_Z(dz)=e^{-1}\sum_{n\ge0}\frac{1}{n!}\delta_n(dz),$ then $c_Z(v)=v\log
    v-v+1$ if $v>0,$ $c_Z(0)=1$ and $c_Z(v)=+\infty$ if $v< 0.$
\end{enumerate}
We have $c_Z(0)=0$ if and only if
$\mathbb{E}Z:=\int_{\RR^d}z\,m_Z(dz)=0.$ More generally,
$c_Z(v)\in [0,+\infty]$ and $c_Z(v)=0$ if
     and only if $v=\mathbb{E}Z.$
We also have
     $$c_{aZ+b}(u)=c_Z\big(a^{-1}(v-b)\big)$$
for all invertible linear operator $a:\RR^d\to\RR^d$ and all
$b\in\RR^d.$
\end{examples}

If $\mathbb{E}Z=0,$ $c_Z$ is quadratic at the origin since
$c_Z(v)= v\cdot\Gamma_Z^{-1}
    v/2+o(|v|^2)$ where $\Gamma_Z$ is the covariance matrix of
    $Z.$ This rules out the usual costs $c(v)=|v|^p$ with $p\not=2.$

Nevertheless, taking $Z$ a real valued variable with  density
    $C\exp(-|z|^p/p)$ with $p\geq 1$ leads to
    $c_Z(v)=|v|^p/p(1+o_{|v|\rightarrow\infty}(1)).$
    The case $p=1$ follows from Example \ref{ex-03}-(3) above. To see that the result still holds with $p>1,$
    compute by means of the Laplace method the principal part as $\zeta$ tends to infinity of $\int_0^\infty
    e^{-z^p/p}e^{\zeta z}\,dz=\sqrt{2\pi(q-1)}\zeta^{1-q/2}e^{\zeta^q/q}(1+o_{\zeta\rightarrow
    +\infty}(1))$ where $1/p+1/q=1.$
    \\
    Of course, we deduce a related $d$-dimensional result considering
    $Z$ with the density
    $C\exp(-|z|_p^p/p)$ where $|z|_p^p=\sum_{i\leq d}|z_i|^p.$ This gives
    $c_Z(v)=|v|_p^p/p(1+o_{|v|\rightarrow\infty}(1)).$

\begin{remark}
Let $R^k$ be defined by \eqref{eq-26} and \eqref{eq-43} where  $Z$
is only allowed to take  isolated values as Examples
\ref{ex-03}-(2) and (4). Suppose that $\mu_0$ has  a discrete
support, then $\pont{R_1}k{\mu_0}$ has
 also a discrete support. It follows that any $P\in\PO$ which is
 absolutely continuous with respect to $\pont{R_1}k{\mu_0}$ is such that $P_1$
 has a discrete support. Now, if you choose a diffuse measure for
 $\mu_1,$ there is no solution to the non-modified minimization
 problem \eqref{eq-35}. We see that it is necessary to introduce a
 sequence $\Seq{\mu_1}k$ of discrete measures such that $\Lim k\mu_1^k=\mu_1$ for the
 sequences of entropy minimization modified problems \eqref{S01k}$_{k\ge1}$ and \eqref{Sk}$_{k\ge1}$ to admit
 solutions.
\end{remark}

\subsection*{Nonlinear transformations}
By means of the contraction principle (Theorem \ref{res-08}), we
can twist the cost functions which have been obtained earlier. We
only present some examples to illustrate this technique.

\subsubsection*{The static case}
Here, we only consider the LD of the final position
$Y^k_1.$ We have just remarked that the cost functions $c_Z$ as
above are necessarily  quadratic at the origin. This  drawback
will be partly overcome by means of  continuous transformations.

We are going to look at an example
$$
    \pont Ykx_1=x+V^k
$$
where $\Seq Vk$ satisfies a LDP which is not given by Cram\'er's
theorem. Let $(Z_j)_{j\geq 1}$ be as above and let $\alpha$ be
    any continuous mapping on $\RR^d.$ Consider
$$
     V^k=\alpha\left(\frac{1}{k}\sum_{1\leq j\leq
    k}Z_j\right).
$$
We obtain $c(v)=\inf\{c_Z(u);u\in \RR^d, \alpha(u)=v\}, v\in\RR^d$
    as a consequence of the contraction principle.
    In particular if $\alpha$ is a continuous injective mapping,
    then
\begin{equation}\label{eq-27}
    c=c_Z\circ \alpha^{-1}.
\end{equation}
For instance, if $Z$ is a standard normal vector as in Example
\ref{ex-03}-(1), we know that the empirical mean of independent
copies of $Z:$ $\frac{1}{k}\sum_{1\leq j\leq
    k}Z_j,$ is a centered normal vector with variance $\mathrm{Id}/k.$ Taking  $\alpha=\alpha_p$
  which is  given for each $p>0$ and $v\in\RR^d$ by
    $\alpha_p(v)=2^{-1/p}|v|^{2/p-1}v,$ leads us to
\begin{equation}\label{eq-29}
    V^k\eqlaw (2k)^{-1/p}|Z|^{2/p-1}Z,
\end{equation}
the equality in law $\eqlaw$ simply means that both sides of the
equality share the same distribution. The mapping $\alpha_p$ has
been chosen to obtain with \eqref{eq-27}:
$$
    c(v):=c_p(v)=|v|^p,\quad v\in\RR^d.
$$
Note that $V^k$ has the same law as $k^{-1/p}Z_p$ where the
density of the law
    of $Z_p$ is $\kappa|z|^{p/2-1}e^{-|z|^p}$ for some normalizing
    constant $\kappa.$

\subsubsection*{The dynamic case}

We now look at an example where
\begin{equation}\label{eq-28}
     \pont Ykx_t=x+V^k_t,\quad 0\le t\le1
\end{equation}
where $\Seq Vk$ satisfies a LDP in $\Omega$ which is not given by
Mogulskii's theorem.

We present examples of dynamics $V^k$ based on the standard
Brownian motion $B=(B_t)_{0\leq t\leq 1}$ in $\RR^d.$ In these
examples, one can restrict the path space to be the space
$\OO=C([0,1],\RR^d)$ equipped with the uniform topology. The item
(1) is already known to us, we recall it for the comfort of the
reader.
\begin{examples}\label{ex-02}\
\begin{enumerate}
\item
    An important example is given by
    \begin{equation*}
    V^k_t =k^{-1/2}B_t,\quad 0\leq t\leq 1.
    \end{equation*}
    Schilder's theorem states that
    $\Seq Vk$ satisfies the LDP in $\Omega$ with the coercive  rate function
\begin{equation*}
    C^0(\omega)=\left\{\begin{array}{ll}
      \int_0^1 |\dot{\omega}_t|^2/2\,dt &  \textrm{if }\omega\in \Oac,\omega_0=0 \\
      +\infty & \textrm{otherwise.} \\
    \end{array}\right.
\end{equation*}
   As in Example \ref{ex-03}-(1), it corresponds to the quadratic cost function
$|v|^2/2,$ but with a different dynamics.
 \item
    More generally, with $p>0,$ we have just seen that
\begin{equation*}
    V^k_t =(2k)^{-1/p}|B_t|^{2/p-1}B_t,\quad 0\leq t\leq 1
\end{equation*}
    corresponds to the power cost function
$ c_p(v)=|v|^p,$  $v\in\RR^d,$ since $V^k_1\eqlaw V^k$ as in
\eqref{eq-29}.
    The associated dynamic cost is given for all $\omega\in\OO$ by
\begin{equation*}
    C^0(\omega)=\left\{\begin{array}{ll}
      p^2/4 \Iii |\omega_t|^{p-2} |\dot{\omega}_t|^2\,dt & \textrm{if }\omega\in \Oac,\omega_0=0 \\
      +\infty & \textrm{otherwise.} \\
    \end{array}\right.
\end{equation*}

\item
    Similarly, with $p>0,$ the dynamics
\begin{equation*}
    V^k_t =(2k)^{-1/p}|B_t/t|^{2/p-1}B_t,\quad 0< t\leq 1
\end{equation*}
    also corresponds to the power cost function
$ c_p(v)=|v|^p,$ $v\in\RR^d,$ since $V^k_1\eqlaw V^k$ as in
\eqref{eq-29}. But, this time the associated dynamic cost is given
for all $\omega\in\OO$ by
\begin{equation*}
    C^0(\omega)=\left\{\begin{array}{ll}
      \frac14 \int_{(0,1]} \1_{\{\omega_t\not=0\}} |\omega_t/t|^{p} \Big|(2-p)\omega_t/|\omega_t|+pt\dot{\omega}_t/|\omega_t|\Big|^2\,dt
      & \textrm{if }\omega\in \Oac, \omega_0=0 \\
      +\infty & \textrm{otherwise.} \\
    \end{array}\right.
\end{equation*}
\end{enumerate}
\end{examples}
Recall that a geodesic path from $x$ to $y$ is some
$\omega\in\Oac$ which solves the minimization problem
\eqref{eq-11}. It is well known that the geodesic paths for Item
(1) are the constant velocity paths $\sigma^{xy},$ see
\eqref{eq-31}. The geodesic paths for Item (2) are still straight
lines but with a time dependent velocity (except for $p=2$). On
the other hand, the geodesic paths for Item (3) are the constant
velocity paths.

\subsection*{Modified random walks on $\RR^d$}
Simple random walks correspond to \eqref{eq-28} with $V^k=W^k$
given by \eqref{eq-43}. We introduce a generalization which is
defined by  \eqref{eq-28} with
\begin{equation*}
    V^k_t=\alpha_t(W^k_t),\quad 0\leq t\leq 1
\end{equation*}
where $\alpha:(t,v)\in [0,1]\times\RR^d \mapsto
\alpha_t(v)\in\RR^d$ is a continuous application such that
$\alpha_0(0)=0$ (remark that $W^k_0=0$ almost surely) and
$\alpha_t$ is injective for all $0<t\leq 1.$
\\
For all $x\in\RR^d$ and all $k\geq 1,$ the random path $\pont
Ykx=x+V^k$ satisfies
$$
\pont Ykx=\Phi(\pont Wkx)
$$
where $\pont Wkx=x+W^k$ and  $\Phi:\OO\to \OO$ is the bicontinuous
injective mapping given for all $\omega\in\OO$ by
$\Phi(\omega)=(\Phi_t(\omega))_{0\leq t\leq 1} $ where
\begin{equation*}
    \Phi_t(\omega)=\omega_0+\alpha_t(\omega_t-\omega_0),
    \quad 0\leq t\leq 1.
\end{equation*}
As for \eqref{eq-27}, the LD rate function of $(\pont
Ykx)_{k\ge1}$ is $C^x=C+\iota_{\{X_0=x\}}$  where
$$
    C=C_Z\circ\Phi^{-1}
$$
and $C_Z$ is given at \eqref{eq-30}. It is easy to see that for
all $\phi\in\OO,$ $\Phi^{-1}(\phi)=(\Phi^{-1}_t(\phi))_{0\le
t\le1}$ where for all $0<t\le1,$
$\Phi^{-1}_t(\phi)=\phi_0+\beta_t(\phi_t-\phi_0)$ with
$\beta_t:=\alpha_t^{-1}.$ Assuming that $\beta$ is differentiable
on $(0,1]\times\RR^d,$ we obtain
\begin{equation*}
    C(\omega)=\left\{\begin{array}{ll}
      \Iii c_Z\big(
    \partial_t\beta_t(\omega_t-\omega_0)+\nabla\beta_t(\omega_t-\omega_0)\cdot\dot\omega_t
      \big)\,dt & \textrm{if }\omega\in \Oac \\
      +\infty & \textrm{otherwise} \\
    \end{array}\right.,\quad \omega\in \OO.
\end{equation*}
For each $x,y\RR^d,$  \eqref{eq-11} admits a unique solution
$\haut\gamma xy$ which is given by the equation
$\Phi^{-1}(\haut\gamma xy)=\pont\sigma x{x+\beta_1(y-x)}$ where
$\haut\sigma xy$ is the constant velocity geodesic, see
\eqref{eq-31}. That is
$$
    \haut\gamma xy_t=x+\alpha_t(t\beta_1(y-x)),
    \quad 0\le t\le1.
$$
The corresponding static cost function $c$ which is specified by
\eqref{eq-10}, i.e.
$$
    c(x,y)=C(\haut\gamma xy),\quad x,y\in\RR^d.
$$ In the case when $\alpha$ doesn't depend on $t,$ we
see that for all $x,y\in\RR^d,$
$$
    c(x,y)=C(\haut\gamma xy)=C_Z(\pont\sigma x{x+\beta(y-x)})=c_Z(\alpha^{-1}(y-x)),
$$
which is \eqref{eq-27}, but the velocity of the geodesic path
    $$\haut{\dot\gamma} xy
    _t=\nabla\alpha\big(t\alpha^{-1}(y-x)\big)\cdot\alpha^{-1}(y-x)$$
is not constant in general.

\section{Proofs of the results of Section \ref{sec-results}}\label{sec-proofs}

The main technical result is Proposition \ref{res-02}.

It will be used at several places that $X_0,X_1:\Omega\to\XX$ are
continuous. This is clear when $\Omega=C(\ii,\XX)$ since it is
furnished with the topology of uniform convergence. In the general
case where $\Omega=D(\ii,\XX)$ is furnished with the Skorokhod
topology, it is known that $X_t$ is not continuous in general.
But, it remains true that $X_0$ and $X_1$ are continuous, due to
the specific form of the metric at the endpoints.

\subsection*{Proof of Proposition \ref{res-02}}
The space $\CO$ is furnished with the supremum norm
$\|f\|=\sup_\Omega |f|,$ $f\in\CO$ and $\CO'$ is its topological
dual space. Let $\MO$, resp.\ $\MOp$ denote the spaces of all
bounded, resp.\ bounded positive, Borel measures on $\Omega.$ Of
course, $\MO\subset\CO'$ with the identification $\langle
f,Q\rangle_{\CO,\CO'}=\IO f\,dQ$ for any $Q\in\MO.$ We write
$\langle f,Q\rangle:=\langle f,Q\rangle_{\CO,\CO'}$ for
simplicity.
\\
Dropping the superscript $k$ for a moment, we have
$(R^x\in\PO;x\in\XX)$ a measurable kernel and $R^{\mu_0}:=\IX
R^x(\cdot)\,\mu_0(dx)$ where $\mu_0\in\PX$ is the initial law.

\begin{lemma}\label{res-16}
For all $ Q\in \CO',$
$$
    H(Q|R^{\mu_0})+\iota_{\{Q\in\PO: Q_0=\mu_0\}}
    =\sup_{f\in\CO}\left\{\langle f,Q\rangle-\IX\log\langle
    e^f,R^x\rangle\,\mu_0(dx)\right\}.
$$
\end{lemma}

This identity should be compared with the well-known variational
representation of the relative entropy
\begin{equation}\label{eq-21}
    H(Q|R)+\iota_{\PO}(Q)=\sup_{f\in\CO}\left\{\langle f,Q\rangle-\log\langle
    e^f,R\rangle\right\},\quad Q\in\MO
\end{equation}
which holds for any reference \emph{probability} measure $R\in\PO$
on any polish space $\Omega.$
\begin{proof}
Denote
$$
    \Theta(f)=\IX\log\langle e^f,R^x\rangle\,\mu_0(dx)\in(-\infty,\infty],
    \quad f\in\CO
$$
Its convex conjugates with respect to the duality $\langle
\CO,\CO'\rangle$ is given for all $Q\in\CO'$ by
    $
    \Theta^*(Q):=\sup_{f\in\CO}\left\{\langle
    f,Q\rangle-\Theta(f)\right\}.
    $
It will be proved at Lemma \ref{res-17} that any $Q\in\CO'$ such
that $\Theta^*(Q)<\infty$ is in $\MOp.$ Let us admit this for a
while, and take $Q\in\MOp$ such that $\Theta^*(Q)<\infty.$ Taking
$f=\phi(X_0)$ with $\phi\in\CX,$ we see that
    $\sup_{\phi\in\CX}\IX \phi\,d(Q_0-\mu_0)\le \Theta^*(Q).$ Hence,
    $\Theta^*(Q)<\infty$ implies that $Q_0=\mu_0.$
This shows us that
    if $\Theta^*(Q)<\infty,$ then $Q$ is a probability measure with $Q_0=\mu_0.$
\\
It remains to prove that for such a $Q\in\PO,$ we have
$\Theta^*(Q)=H(Q|R^{\mu_0}).$ Since $\OO$ is a polish space, any
$Q\in\PO$ such that $Q_0=\mu_0$ disintegrates as
$$Q(\cdot)=\IX Q^x(\cdot)\,\mu_0(dx)$$ where $(Q^x;x\in\XX)$ is a
measurable kernel of probability measures. We see that
    $$\Theta^*(Q)=\sup_{f\in\CO}\IX[\langle f,Q^x\rangle-\log\langle e^f,R^x\rangle]\,\mu_0(dx).$$
We obtain
\begin{eqnarray*}
  \Theta^*(Q)
  &\le& \IX \sup_{f\in\CO}[\langle f,Q^x\rangle-\log\langle e^f,R^x\rangle]\,\mu_0(dx) \\
  &\overset{\checkmark}=& \IX H(Q^x|R^x)\,\mu_0(dx) \\
  &=&  H(Q|R^{\mu_0})
\end{eqnarray*}
where \eqref{eq-21} is used at the marked equality and last
equality follows from the tensorization property \eqref{eq-20}.
Note that $x\mapsto H(Q^x|R^x)$ is measurable. Indeed,
$(Q,R)\mapsto H(Q|R)$ is \lsc\ being the supremum of continuous
functions, see \eqref{eq-21}. Hence, it is Borel measurable. On
the other hand, $x\mapsto R^x$ and $x\mapsto Q^x$ are also
measurable, being the disintegration kernels of   Borel measures
on a polish space.

Let us prove the converse inequality. By Jensen's inequality:
    $
    \IX \log\langle e^f,R^x\rangle\,\mu_0(dx)
    \le \log\IX \langle e^f,R^x\rangle\,\mu_0(dx)
    =\log \langle e^f,R^{\mu_0}\rangle,
    $
so that
$$
 \Theta^*(Q)\ge\sup_{f\in\CO}\left\{\IO f\,dQ-\log \IO e^f\,dR^{\mu_0}\right\}
 =H(Q|R^{\mu_0})
$$
where the equality is \eqref{eq-21} again. This completes the
proof of the lemma.
\end{proof}

During the proof of Lemma \ref{res-16}, we used a result which is
stated at Lemma \ref{res-17}. Denote also
$$
    \Lambda(f):=\IX\sup_\Omega\{f-C^x\}\,\mu_0(dx)=\IX\sup_{\Omega^x}\{f-C\}\,\mu_0(dx),
    \quad f\in\CO
$$
where $\Omega^x:=\{X_0=x\}\subset\Omega.$ It will appear later
that the function $\Lambda$ is the convex conjugate of the
$\Gamma$-limit $\mathcal{C}.$ Its convex conjugate with respect to
the duality $\langle \CO,\CO'\rangle$ is given for all $Q\in\CO'$
by
    $
     \Lambda^*(Q):=\sup_{f\in\CO}\left\{\langle
    f,Q\rangle-\Lambda(f)\right\}.
    $
\begin{lemma}\label{res-17}\
\begin{enumerate}
    \item  $\{\Theta^*<\infty\}\subset \MOp;$
    \item $\{\Lambda^*<\infty\}\subset \MOp.$
\end{enumerate}
\end{lemma}

\proof For a positive element $Q\in\CO'$ to be in $\MO,$ it
necessary and sufficient that it is $\sigma$-additive. That is,
for all \emph{decreasing} sequence $\seq fn$ in $\CO$ such that
$\Lim n f_n=0$ \emph{pointwise}, we have $\Lim n \langle
f_n,Q\rangle=0.$

\Boulette{(1)} Let us prove that $\{\Theta^*<\infty\}\subset
\MOp.$
\\
Let us show that $Q\geq 0$ if $\Theta^*(Q)<\infty.$ Let $f\in \CO$
be such that $f\geq 0.$ As $\Theta( af)\leq 0$ for all $a\leq 0,$
\begin{eqnarray*}
   \Theta^*(Q)
  &\geq& \sup_{a\leq 0}\{a\langle f,Q\rangle-\Theta(af)\} \\
  &\geq& \sup_{a\leq 0}\{a\langle f,Q\rangle\}\\
  &=& \left\{%
\begin{array}{ll}
    0, & \hbox{if } \langle f,Q\rangle\geq 0\\
    +\infty, & \hbox{otherwise.} \\
\end{array}%
\right.
\end{eqnarray*}
Therefore, if $\Theta^*(Q)<\infty,$ $\langle f,Q\rangle\geq 0$ for
all $f\geq 0,$ which is the desired result.

Let us take a decreasing sequence  $\seq fn$ in $\CO$ which
converges pointwise to zero. By the dominated convergence theorem,
we have
$$
    \Lim n \Theta(af_n)=0,\quad \forall a\ge0.
$$
It follows that for all $Q\in \CO',$
\begin{eqnarray*}
  \Theta^*(Q)
  &\geq& \sup_{a\geq 0}\limsup_{n\rightarrow\infty}\{a\langle f_n,Q\rangle-\Theta(af_n)\} \\
   &\geq& \sup_{a\geq 0}\left(\limsup_{n\rightarrow\infty} a\langle f_n,Q\rangle-\lim_{n\rightarrow\infty}\Theta(af_n)\right) \\
   &=& \sup_{a\geq 0}a\limsup_{n\rightarrow\infty} \langle f_n,Q\rangle \\
  &=& \left\{%
\begin{array}{ll}
    0 & \hbox{if }\limsup_{n\rightarrow\infty} \langle f_n,Q\rangle\leq 0 \\
    +\infty & \hbox{otherwise.} \\
\end{array}%
\right.
\end{eqnarray*}
Therefore, if $\Theta^*(Q)<\infty,$ we have
$\limsup_{n\rightarrow\infty} \langle f_n,Q\rangle\leq 0.$ Since
we have just seen that $Q\geq 0,$ we have the desired result.

\Boulette{(2)} Let us prove that $\{\Lambda^*<\infty\}\subset
\MOp.$
\\
Let us show that $Q\geq 0$ if $\Lambda^*(Q)<\infty.$ Let $f\in
\CO$ be such that $f\geq 0.$ As $\inf C=0,$ $\Lambda( af)\leq 0$
for all $a\leq 0,$ and we conclude as at item (1).

Let us take a decreasing sequence  $\seq fn$ in $\CO$ which
converges pointwise to zero. By Lemma \ref{res-18} below, for all
$x\in\XX,$ $\big(\sup_\Omega \{f_n-C^x\}\big)_{n\geq 1}$ is a
decreasing sequence and $\Lim n\sup_\Omega \{f_n-C^x\}=0.$ As
$\left|\sup_\Omega \{f_n-C^x\}\right|\leq \sup_\Omega|f_1|<\infty$
for all $n$ and $x,$ we can apply the dominated convergence
theorem to obtain that $ \Lim n\Lambda(af_n)=0, $ for all $a\geq
0$ and we conclude as at item (1).

Finally, one must be careful with the measurability of
$x\in\XX\mapsto
u_n(x):=\inf_\Omega\{C^x-f_n\}=-\sup_\Omega\{f_n-C^x\}\in\RR.$
Since $\Omega$ and $\XX$ are assumed to be polish, we can apply a
general result by Beiglb\"ock and Schachermayer \cite[Lemmas 3.7,
3.8]{BeSc08} which tells us that for each $n\ge1$ and each Borel
probability measure $\mu$ on $\XX,$ there exists a Borel
measurable function $\tilde{u}_n$ on $\XX$ such that
$\tilde{u}_n\le u_n$ and $\tilde{u}_n(x)= u_n(x)$ for $\mu$-a.e.
$x\in\XX.$
\endproof

During the proof of the previous lemma we have invoked the
following result.

\begin{lemma}\label{res-18}
Let $J$ be a coercive $[0,\infty]$-valued function on $\Omega$ and
$(f_n)_{n\geq 1}$ a decreasing sequence of continuous bounded
functions on $\Omega$ which converges pointwise to some bounded
upper semicontinuous function $f.$ Then,
$\left(\sup_\Omega\{f_n-J\}\right)_{n\geq 1}$ is a decreasing
sequence and
$$
\Lim n \sup_\Omega\{f_n-J\}=\sup_\Omega\{f-J\}.
$$
\end{lemma}

 \proof
Changing sign and denoting $g_n= J-f_n,$ $g=J-f,$ we want to prove
that $\Lim n\inf_\Omega g_n=\inf_\Omega g.$

We see that $(g_n)_{n\geq 1}$ is an increasing sequence of \lsc\
functions. It follows by the Proposition 5.4 of \cite{DalMaso}
that it is a $\Gamma$-convergent sequence and
\begin{equation}\label{eq-23a}
    \Glim n g_n=\Lim n g_n=g.
\end{equation}
Let us admit for a while that there exists some compact set $K$
which satisfies
\begin{equation}\label{eq-23b}
     \inf_\Omega g_n= \inf_Kg_n
\end{equation}
for all $n.$ This and the convergence (\ref{eq-23a}) allow to
apply Theorem 7.4 of \cite{DalMaso} to obtain $\Lim n \inf_\Omega
g_n= \inf_\Omega \Glim n g_n=\inf_\Omega g$ which is the desired
result.

\par\medskip
It remains to check that (\ref{eq-23b}) is true. Let
$\omega_*\in\Omega$ be such that $J(\omega_*)<\infty$ (if $J\equiv
+\infty,$ there is nothing to prove). Then, $\inf_\Omega g_n\leq
g_n(\omega_*)=J(\omega_*)-f_n(\omega_*)\leq
J(\omega_*)-f(\omega_*) \leq J(\omega_*)-\inf_\Omega f.$ On the
other hand, for all  $n,$ $f_n\leq f_1 \leq A:= \sup f_1.$ Let
$B:= A+1+J(\omega_*)-\inf_\Omega f.$ For all $\omega$ such that
$J(\omega)>B,$ we have $g_n(\omega)>B-\sup_\Omega f_n\geq B-A\geq
J(\omega_*)-\inf_\Omega f+1.$ We have just seen that for all $n,$
\begin{equation*}
    \inf_\Omega g_n\leq J(\omega_*)-\inf_\Omega f \qquad \textrm{and}
    \qquad
    \inf_{\omega; J(\omega)>B} g_n(\omega)\geq J(\omega_*)-\inf_\Omega
    f+1.
\end{equation*}
This proves (\ref{eq-23b}) with the compact level set $K=\{J\leq
B\}$ and completes the proof of the lemma.
\endproof

Recall that for all $Q\in\MO,$ $ \mathcal{C}^{k,\mu_0}(Q)=\frac1k
H(Q|\pont R{k}{\mu_0})+\iota_{\{Q_0=\mu_0\}}.$ With Lemma
\ref{res-16}, we see that
\begin{equation}\label{eq-19}
      \mathcal{C}^{k,\mu_0}(Q)=\Lambda_k^*(Q),\quad Q\in\CO'
\end{equation}
where $\Lambda_k^*$ is the convex conjugate of
$$
    \Lambda_k(f)=\IX\frac1k\log\langle e^{kf},\pont Rkx\rangle\,
    \mu_0(dx),
    \quad f\in\CO
$$
with respect to the duality $\langle\CO,\CO'\rangle.$ The keystone
of the proof of Proposition \ref{res-02} is the following
consequence of the Laplace-Varadhan principle.

\begin{lemma}\label{res-19}
Under the assumptions of Proposition \ref{res-02}, for all
$f\in\CO,$ we have
\begin{enumerate}
    \item $\Lim k\Lambda_k(f)=\Lambda(f);$
    \item $\sup_{k\ge1} |\Lambda_k(f)|\leq \|f\|,\quad |\Lambda(f)|\leq \|f\|:= \sup_\Omega|f|.$
\end{enumerate}
The functions $\Lambda_k$ and $\Lambda$ are convex.
\end{lemma}

\proof Our assumptions allow us to apply the Laplace-Varadhan
principle, see  Theorem \ref{res-15}. It tells us that for each
$x\in\XX,$
$$
    \Lim k\frac{1}{k} \log \langle  e^{kf},\pont Rkx\rangle= \sup_\Omega\{f-C^x\}.
$$
On the other hand, it is clear that for each $k\ge1,$
$|\frac{1}{k} \log \langle  e^{kf},\pont Rkx\rangle|\le\|f\|.$
Passing to the limit, we also get
$|\sup_\Omega\{f-C^x\}|\le\|f\|.$ Now by the Lebesgue dominated
convergence theorem, we obtain the statements (1) and (2).
\\
Note that $x\mapsto \sup_\Omega\{f-C^x\}$ is measurable as a
pointwise limit of measurable functions.
\\
It is standard to prove with H\"older's inequality that $f\mapsto
\frac{1}{k} \log \langle  e^{kf},\pont Rkx\rangle$ is convex. It
follows that $\Lambda_k$ and $\Lambda$ are also convex.
\endproof

We are in position to apply Corollary \ref{res-13}. Let us equip
$\CO'$ with the $*$-weak topology $\sigma(\CO',\CO)$. By Corollary
\ref{res-13}, we have
\begin{equation}\label{eq-18}
    \Glim k \Lambda_k^*=\Lambda^*
\end{equation}
where
$$
    \Lambda^*(Q)=\sup_{f\in\CO} \left\{\langle f,Q\rangle_{\CO,\CO'}-\IX\sup_{\Omega}\{f-C^x\}\,\mu_0(dx)\right\},
    \quad Q\in \CO'.
$$
This limit still holds in $\MO\subset\CO',$ by Lemma \ref{res-17}.

Because of \eqref{eq-19}, \eqref{eq-18} and Lemma \ref{res-17}, to
complete the proof of Proposition \ref{res-02}, it remains to
prove the subsequent lemma.
\begin{lemma}
Let $C$ be a \lsc\ $[0,\infty]$-valued function on the polish
space $\Omega$. Denote $C^x=C+\iota_{\{\theta=x\}}$ for each
$x\in\XX,$ where $\theta:\Omega\to\XX$ is a continuous application
with its values in polish space $\XX$. Take $\mu\in\PX$ and
suppose that
$$\inf_\Omega C^x=0$$ for $\mu$-almost every $x\in\XX.$ Then, we
have
\begin{multline}\label{eq-22}
   \qquad\sup_{f\in\CO} \left\{\langle f,Q\rangle-\IX\sup_{\Omega}\{f-C^x\}\,\mu(dx)\right\}
    \\ =\IO C\,dQ+\iota_{\{Q\in\PO:\theta\pf Q=\mu\}},
    \quad Q\in \MO.\qquad
\end{multline}
\end{lemma}
Note that since $C\ge0$ and $C$ is measurable, the integral $\IO
C\,dP$ makes sense in $[0,\infty]$ for any $P\in\PO.$
\\
As the function $C$ of Proposition \ref{res-02} is such that $C^x$
is a LD rate function for all $x\in\XX,$ it satisfies the
assumption $\inf_\Omega C^x=0$ for $\mu$-almost every $x\in\XX.$

\proof  Let us first check that if $Q\in\MO$ satisfies
$\Lambda^*(Q)<\infty,$ then $Q\in\PO$ and $\theta\pf Q=\mu\in\PX.$
We already know by Lemma \ref{res-17} that $Q\in\MOp.$  Choosing
$f=\phi\circ\theta$ with $\phi\in\CX,$ since $\inf_\Omega C^x=0,$
we see that
    $\sup_{\Omega}\{\phi\circ\theta-C^x\}=\phi(x).$ Hence, $\sup_{\phi\in\CX}\IX
    \phi\,d(\theta\pf Q-\mu)\le \Lambda^*(Q)<\infty$
    which implies that $\theta\pf Q=\mu.$ This proves the desired result.

It remains to prove the equality for a fixed $P\in\PO$ which
satisfies $\theta\pf P=\mu.$ Because $\Omega$ and $\XX$ are polish
spaces, we know that $P$ disintegrates as follows:
    $P(\cdot)=\IX P^x(\cdot)\,\mu(dx),$ with $x\in\XX\mapsto P^x(\cdot):=P(\cdot\mid
    \theta=x)\in\PO$ Borel measurable.
For any $f\in\CO,$
\begin{eqnarray*}
  \langle f,P\rangle-\IX\sup_\Omega\{f-C^x\}\,\mu(dx)
  &=& \IX [\langle f,P^x\rangle-\sup_\Omega\{f-C^x\}]\,\mu(dx)\\
  &=& \IX [\langle C^x,P^x\rangle+\langle f-C^x-\sup_\Omega\{f-C^x\},P^x\rangle]\,\mu(dx) \\
  &\le& \IX \langle C^x,P^x\rangle\,\mu(dx) \\
  &=& \IO C\,dP.
\end{eqnarray*}
 Optimizing, we obtain
\begin{equation*}
     \sup_{f\in\CO}\left\{\langle f,P\rangle-\IX\sup_\Omega\{f-C^x\}\,\mu(dx)\right\}
     \le\IO C\,dP.
\end{equation*}
If $C$ is in $\CO,$ the case of equality is obtained with $f=C,$
$P$-a.e.\ and in this situation we see that the identity
\eqref{eq-22} is valid. This will be invoked very soon.

In the general case, $C$ is only assumed to be \lsc. By means of
the Moreau-Yosida approximation procedure which is implementable
since $\Omega$ is a metric space, one can build an increasing
sequence $\seq Cn$ of functions in $\CO$ which converges pointwise
to $C.$ Therefore,
\begin{eqnarray*}
   &&\sup_{f\in\CO}\left\{\langle f,P\rangle-\IX\sup_\Omega\{f-C^x\}\,\mu(dx)\right\}\\
   &\le& \IO C\,dP \\
  &\overset{(\textrm{i})}=& \sup_{n\ge1}  \IO C_n\,dP\\
  &\overset{(\textrm{ii})}=& \sup_{n\ge1} \sup_{f\in\CO}\left\{\langle f,P\rangle-\IX\sup_\Omega\{f-C_n^x\}\,\mu(dx)\right\}\\
  &=&  \sup_{f\in\CO}\left\{\langle f,P\rangle+\sup_{n\ge1}\IX\inf_\Omega\{C_n^x-f\}\,\mu(dx)\right\}\\
  &\overset{(\textrm{iii})}\le&\sup_{f\in\CO}\left\{\langle f,P\rangle+\IX\inf_\Omega\{C^x-f\}\,\mu(dx)\right\}\\
  &=& \sup_{f\in\CO}\left\{\langle f,P\rangle-\IX\sup_\Omega\{f-C^x\}\,\mu(dx)\right\},
\end{eqnarray*}
which proves the desired identity \eqref{eq-22}.
\\
Equality (i) follows from the monotone convergence theorem. Since
$C_n$ stands in $\CO,$ equality (ii) is valid  (this has been
proved a few lines earlier) and the inequality (iii) is a direct
consequence of $C_n\le C$ for all $n\ge1.$ Note that
$x\in\XX\mapsto \inf_\Omega\{C_n^x-f\}\in\RR$ is upper
semicontinuous and it is a fortiori Borel measurable.
\endproof

\subsection*{Proofs of the remaining results}
The keystone of the proofs of the remaining results is Proposition
\ref{res-02}.
\subsubsection*{Proposition \ref{res-05}}
Proposition \ref{res-05} is a particular case of Proposition
\ref{res-02}. Indeed, choosing $\Omega=\XXX$ which can be
interpreted as the space of all $\XX$-valued paths on the
two-point time interval $\{0,1\},$ and taking
$C(\omega)=c(\omega_0,\omega_1)$ where $c$ is assumed to be \lsc,
with $\omega=(x,y)$  we see that
$C^x(x',y)=c(x,y)+\iota_{\{x'=x\}}$ for all $x,x',y\in\XX.$ The
assumption that $c(x,\cdot)$ is coercive on $\XX$ is equivalent to
the coerciveness of $C^x$ on $\XXX.$

\subsubsection*{Corollary \ref{res-03} and Theorem \ref{res-04}}
With Proposition \ref{res-02} in hand, Corollary \ref{res-03} and
Theorem \ref{res-04} are immediate consequences of Theorem
\ref{res-10} and of the equi-coerciveness with respect to the
$*$-weak topology $\sigma(\PO,\CO)$ of
$\{\mathcal{C},\mathcal{C}^k;k\ge 1\}.$ This equi-coerciveness
follows from Corollary \ref{res-13} and Lemma \ref{res-19}. The
uniqueness of the solution to \eqref{Sk} follows from the strict
convexity of the relative entropy.

\subsubsection*{Corollary \ref{res-06} and Theorem \ref{res-07}}
Similarly, once we have Proposition \ref{res-05} in hand,
Corollary \ref{res-06} and Theorem \ref{res-07} are immediate
consequences of Theorem \ref{res-10} and of the equi-coerciveness
with respect to the $*$-weak topology $\sigma(\PX,\CX)$ of
$\{\bas{\mathcal{C}}01,\bas{\mathcal{C}^k}01;k\ge 1\}.$ This
equi-coerciveness follows from the fact that the set of all
probability measures $\pi\in\PXX$ such that $\pi_0=\mu_0$ and
$\pi_1\in\{\mu_1,\mu_1^k;k\ge1\}$ is relatively compact since
$\Lim k\mu_1^k=\mu_1$; a consequence of Prokhorov's theorem in a
polish space.
\\
Again, the uniqueness of the solution to \eqref{S01k} follows from
the strict convexity of the relative entropy.
\\
Note that, when $C$ and $c$ are linked by \eqref{eq-10}, one can
also derive the equi-coerciveness of
$\{\bas{\mathcal{C}}01,\bas{\mathcal{C}^k}01;k\ge 1\}$ from the
equi-coerciveness  of $\{\mathcal{C},\mathcal{C}^k;k\ge 1\},$ as
in the proof of Theorem \ref{res-10}.

\subsubsection*{Theorem \ref{res-09}}
The proof of Theorem \ref{res-09} relies upon the subsequent
lemma.

\begin{lemma}\label{res-20}
Under the assumptions of Proposition \ref{res-02}, the function
$c$ defined by \eqref{eq-10} is \lsc\ and
$$\inf\left\{\IO C\,dP; P\in\PO,\bas P01=\pi\right\}
     =\IX c\,d\pi\in[0,\infty],$$ for all $\pi\in\PXX.$
\end{lemma}

\proof Let us define the function
$$
    \Psi(\pi):=\inf\left\{\IO C\,dP;P\in\PO: \bas P01=\pi\right\},
    \quad \pi\in\PXX.
$$
As $C$ is assumed to be \lsc\ on $\OO,$ $\Psi$ satisfies the
Kantorovich type  dual equality:
\begin{equation}\label{eq-32}
    \Psi(\pi)=\sup_{f\in\mathcal{F}}\IXX f\,d\pi,
    \quad\pi\in\PXX
\end{equation}
where $\mathcal{F}:=\{f\in\CXX; f(X_0,X_1)\le C\}.$ For a proof of
\eqref{eq-32}, one can rewrite mutatis mutandis the proof of the
Kantorovich dual equality. See for instance \cite[Thm 3.2]{Leo07b}
and note that this result takes into account cost functions which
may take infinite values as in the present case.
\\
This shows that $\Psi$ is a \lsc\ function on $\PXX,$ being the
supremum of continuous functions. Define the function
$$
    \psi(x,y):=\Psi(\delta_{(x,y)}),
    \quad x,y\in\XX.
$$
We deduce immediately from the lower semicontinuity of $\Psi$ that
$\psi$ is \lsc\ on $\XXX$. Hence it is Borel measurable. Since it
is $[0,\infty]$-valued, the integral $\IX \psi\,d\pi$ is
meaningful for all $\pi\in\PXX.$ We are going to prove that
\begin{equation}\label{eq-33}
    \Psi(\pi)=\IXX \psi\,d\pi,\quad\pi\in\PXX.
\end{equation}
For any $\pi\in\PXX,$ we obtain
\begin{eqnarray*}
  \Psi(\pi)
  &=& \inf\left\{\IXX\left(\IO C\,d\haut Pxy\right)\,\pi(dxdy);P\in\PO\right\} \\
  &\ge& \IXX  \inf\left\{\IO C\,dP;P\in\PO:\bas P01=\delta_{(x,y)}\right\}\,\pi(dxdy) \\
  &=& \IXX \psi\,d\pi.
\end{eqnarray*}
Let us show the converse inequality. With \eqref{eq-32}, we see
that for each $f\in\mathcal{F}$ and all $(x,y)\in\XXX,$
$\psi(x,y)=\Psi(\delta_{(x,y)})\ge \IXX
f\,d\delta_{(x,y)}=f(x,y).$ That is $f\le\psi,$ for all
$f\in\mathcal{F}.$ Therefore, $
\Psi(\pi)=\sup_{f\in\mathcal{F}}\IXX f\,d\pi\le\IXX\psi\,d\pi,$
completing the proof of \eqref{eq-33}.

It remains to establish that $\psi=c.$ With \eqref{eq-32}, we get
$\psi=\sup \mathcal{F}.$ But it is clear that $f\in\mathcal{F}$ if
and only if for all $x,y\in\XX,$ $f(x,y)\le
\inf\{C(\omega);\omega\in\OO:\omega_0=x,\omega_1=y\}:=c(x,y).$
Hence, $\psi$ is the upper envelope of the set of all functions
$f\in\CXX$ such that $f\le c.$ In other words $\psi$ is the \lsc\
envelope $\ls c$ of $c.$ Finally, for all $x,y\in\XX,$
    $\ls c(x,y)=\psi(x,y)=\inf\left\{\IO C\,d\haut Pxy;P\in\PO\right\}
    \ge c(x,y)\ge \ls c(x,y).$ This implies the desired result: $\psi=\ls c=c.$
\endproof

With this result at hand, let us prove Theorem \ref{res-09}. It is
assumed that for any $x\in\XX,$ $(\pont Rkx)_{k\ge1}$ satisfies
the LDP with scale $k$ and rate function $C^x.$ We have $\pont
\rho kx=(X_1)\pf\pont Rkx.$ Taking the continuous image
$X_1:\Omega\to\XX,$ by means of the contraction principle, see
Theorem \ref{res-08} at the Appendix, we obtain that for any
$x\in\XX,$ $(\pont \rho kx)_{k\ge1}$ satisfies the LDP with scale
$k$ and rate function
$$
    y\in\XX\mapsto
    \inf\{C^x(\omega);\omega\in\Omega:\omega_1=y\}=c(x,y)\in[0,\infty].
$$
    \Boulette{(1)} The first assertion of Theorem \ref{res-09} follows
from the lower semicontinuity of $c$ which was obtained at Lemma
\ref{res-20}. Indeed, this shows that the assumptions of
Proposition \ref{res-05} are fulfilled. The identity
$\inf\eqref{MKdyn}=\inf\eqref{MK}$ is a direct consequence of
Lemma  \ref{res-20}.
    \Boulette{(2)}The second
assertion follows from $\inf\eqref{MKdyn}=\inf\eqref{MK},$ the
convergence of the minimal values which was obtained at item (1)
together with the strict convexity (for the uniqueness) and the
coerciveness (for the existence) of the relative entropy. The
relation between $\Ph^k$ and $\pih^k$ is \eqref{eq-06}.

    \Boulette{(3)} Let us first show that
$P\mapsto\langle C,P\rangle+\iota_{\{P_0=\mu_0\}}$ is coercive on
$\PO.$ By \eqref{eq-22} and the proof of Corollary \ref{res-13},
we see that its sublevel sets  are relatively compact. Since $C$
is \lsc, it is also \lsc. Therefore, it is coercive and so is
$P\mapsto\langle C,P\rangle+\iota_{\{P_0=\mu_0,P_1=\mu_1\}}$. In
particular, if $\inf\eqref{MKdyn}<\infty,$ the set of minimizers
of \eqref{MKdyn} is a nonempty convex compact subset of $\PO.$
\\
Let $\Ph$ be such a minimizer. It disintegrates as
$\Ph(\cdot)=\IXX \haut\Ph xy(\cdot)\,\bas\Ph01(dxdy)$ and with
Lemma \ref{res-20}, we see that $\bas\Ph01:=\pih$ is a solution to
\eqref{MK}. Moreover,
    $
    \IXX c\,d\pih=\psi(\pih)=\IO C\,d\Ph=\IXX
    \left(\IO C\,d\haut\Ph xy\right)\,\pih(dxdy)
    $
and $\IO C\,d\haut\Ph xy\ge c(x,y)$ for $\pih$-a.e.\ $(x,y).$
Hence, $\IO C\,d\haut\Ph xy= c(x,y)$ for $\pih$-a.e.\ $(x,y).$
This means that for $\pih$-a.e.\ $(x,y),$ $\haut\Ph xy(\haut\Gamma
xy)=1$ where $\haut\Gamma xy:=\{\omega\in\OO;
\omega_0=x,\omega_1=y,C(\omega)=c(x,y)\}$ is the set of all
geodesic paths from $x$ to $y.$ Remark that $\haut\Gamma xy$ is a
compact subset of $\OO$ which is nonempty as soon as
$c(x,y)<\infty.$ In particular, it is a Borel measurable subset.
Following the cases of equality, it is clear that if, conversely
$P\in\PO$ satisfies $\haut Pxy(\haut\Gamma xy)=1$ for $\bas
P01$-a.e.\ $(x,y),$ then $P$ minimizes $Q\mapsto\IO C\,dQ$ subject
 to $\bas Q01=\bas P01.$ This completes the proof of the theorem.

\section{$\Gamma$-convergence of convex functions on a weakly compact space}\label{sec-Gamma1}

A typical result about the $\Gamma$-convergence of a sequence of
convex functions $\seq fk$ is: If the sequence of the convex
conjugates $\seq{f^*}k$ converges in some sense, then $\seq fk$
$\Gamma$-converges. Known results of this type are usually stated
in separable reflexive Banach spaces. For instance Corollary 3.13
of H. Attouch's monograph \cite{Attouch84} is

\begin{theorem}\label{res-Attouch}
Let $X$ be a separable reflexive Banach space and $\seq fk$ a
sequence of closed convex functions from $X$ into
$(-\infty,+\infty]$ satisfying the equi-coerciveness assumption: $
f_k(x)\geq\alpha(\|x\|)$ for all $x\in X$ and $k\geq 1$ with
$\lim_{r\rightarrow +\infty} \alpha(r)/r=+\infty. $ Then, the
following statements are equivalent
\begin{enumerate}
    \item   $f=\mathrm{seq }X_w\textrm{-}\Gamma\textrm{-}\Lim k f_k$
    \item $f^*=X_s^*\textrm{-}\Gamma\textrm{-}\lim_{n\rightarrow
    \infty}f_k^*$
    \item $\forall y\in X^*,$ $f^*(y)=\Lim kf_k^*(y)$
\end{enumerate}
where $X^*$ is the dual space of $X,$ $\mathrm{seq }X_w$ refers to
the weak sequential convergence in $X$ and $X_s^*$ to the strong
convergence in $X^*.$
\end{theorem}

Going beyond the reflexivity assumption is not so easy, as can be
seen in Beer's monograph \cite{Beer93}.

In some applications in probability, the reflexive Banach space
setting is not as natural as it is for the usual applications of
variational convergence to PDEs. For instance when dealing with
random measures on $\XX,$ the narrow topology $\sigma(\PX,
C_b(\XX))$ doesn't fit the above framework since $C_b(\XX)$
endowed with the uniform topology may not be separable (unless
$\XX$ is compact) and is not reflexive.

 The next result is an analogue of
Theorem \ref{res-Attouch} which agrees with applications for
random probability measures. Since we didn't find it in the
literature, we give its detailed proof.

Let $X$ and $Y$ be two vector spaces in separating duality. The
space $X$ is furnished with the weak topology $\sigma(X,Y).$

We denote $\iota_C$ the indicator function of the subset $C$ of
$X$ which is defined by $\iota_C(x)=0$ if $x$ belongs to $C$ and
$\iota_C(x)=+\infty$ otherwise. Its convex conjugate is the
support function of $C:$ $\iota_C^*(y)=\sup_{x\in C}\langle
x,y\rangle,$ $y\in Y.$

\begin{theorem}\label{res-41}
Let $\seq gk$ be a sequence of functions on $Y$ such that
\begin{itemize}
    \item[(a)] for all $k,$ $g_k$ is a real-valued convex function on
$Y,$
    \item[(b)] $\seq gk$ converges pointwise to $g:=\Lim k g_k,$
    \item[(c)] $g$ is real-valued and
    \item[(d)] in restriction to any finite dimensional vector subspace $Z$ of $Y,$
    $\seq gk$ $\Gamma$-converges to $g,$ i.e. $\Glim k
    (g_k+\iota_Z)=g+\iota_Z,$ where $\iota_Z$ is the indicator
    function of $Z.$
\end{itemize}
Denote the convex conjugates on $X:$ $f_k=g^*_k$ and $f=g^*.$

If in addition,
\begin{itemize}
    \item[(e)] there exists a $\sigma(X,Y)$-compact set $K\subset X$ such that
    $\dom f_k\subset K$ for all $k\geq 1$ and $\dom f\subset K$
\end{itemize}
then, $\seq fk$ $\Gamma$-converges to $f$ with respect to
$\sigma(X,Y).$
\end{theorem}

The proof of this theorem is postponed after the two preliminary
Lemmas \ref{res-A1} and \ref{res-43}.

\begin{remark}\label{A09}
By (\cite{DalMaso}, Proposition 5.12), under the assumption (a),
assumption (d) is implied by:
\begin{itemize}
    \item[(d')] in restriction to any finite dimensional vector subspace $Z$ of $Y,$
    $\seq gk$ is equibounded, i.e. for all $y_o\in Z,$ there exists $\delta >0$
    such that
\begin{equation*}
    \sup_{k\geq 1}\sup\{|g_k(y)|; y\in Z, |y-y_o|\leq
    \delta\}<\infty.
\end{equation*}
\end{itemize}
\end{remark}

A useful consequence of Theorem \ref{res-41} is

\begin{corollary}\label{res-13}
Let $(Y,\|\cdot\|)$ be a normed space and $X$ its topological dual
space. Let $\seq gk$ be a sequence of functions on $Y$ such that
\begin{itemize}
    \item[(a)] for all $k,$ $g_k$ is a real-valued convex function on $Y,$
    \item[(b)] $\seq gk$ converges pointwise to $g:=\Lim k g_k$ and
    \item[(d'')]  there exists $c>0$ such that $|g_k(y)|\leq
    c(1+\|y\|)$ for all $y\in Y$ and $k\geq 1.$
\end{itemize}
Then, $\seq fk$ $\Gamma$-converges to $f$ with respect to
$\sigma(X,Y)$ where $f_k=g^*_k$ and $f=g^*.$
\\
Moreover, there exists a $\sigma(X,Y)$-compact set $K\subset X$
such that
    $\dom f_k\subset K$ for all $k\geq 1$ and $\dom f\subset K.$
\end{corollary}

 \proof
Under (b), (d'') implies (c). As (d'') implies (d'), we have (d)
by Remark \ref{A09}. Finally, (d'') implies (e) with $K=\{x\in X;
\|x\|_*\leq c\}$ where $\|x\|_*=\sup_{y, \|y\|\leq 1}\langle
x,y\rangle$ is the dual norm on $X.$ Indeed, suppose that for all
$y\in Y,$ $g(y)\leq c+c\|y\|$ and take $x\in X$ such that
$g^*(x)<+\infty.$ As for all $y,$ $\xy \leq g(y)+ g^*(x),$ we get
$|\xy|/\|y\|\leq (g^*(x)+c)/\|y\|+c.$ Letting $\|y\|$ tend to
infinity gives $\|x\|_*\leq c$ which is the announced result.
\\
The conclusion follows from Theorem \ref{res-41}.
\endproof

\begin{lemma}\label{res-A1}
Let $f:X\rightarrow (-\infty,+\infty]$ be a \lsc\ convex function
such that $\dom f$ is included in a compact set. Let $V$ be a
closed convex subset of $X.$

Then, if $V$ satisfies
\begin{equation}\label{A01}
    V\cap \dom f\not=\emptyset \quad \textrm{or}\quad V\cap
    \cl\dom f=\emptyset,
\end{equation}
we have
\begin{equation}\label{A02}
    \infV f(x)=-\infY (f^*(y)+\iota^*_V(-y))\in(-\infty,\infty]
\end{equation}
and if $V$ doesn't satisfy (\ref{A01}), we have
\begin{equation}\label{A03}
    \inf_{x\in W}f(x)=-\infY(f^*(y)+\iota^*_{W}(-y))=+\infty
\end{equation}
for all closed convex set $W$ such that $W\subset\inter V.$
\end{lemma}

 \proof
The proof is divided in two parts. We first consider the case
where $V\cap \dom f\not=\emptyset,$ then the case where $V\cap
\cl\dom f=\emptyset.$

$\bullet$ \textit{The case where $V\cap \dom f\not=\emptyset.$} As
$V$ is a nonempty closed convex set, its indicator function
$\iota_V$ is a closed convex function so that its biconjugate
satisfies $\iota_V^{**}=\iota_V,$ i.e. $\iota_V(x)=\supY
\{\xy-\iota_V^*(y)\}$ for all $x\in X.$ Consequently,
\begin{equation*}
    \infV f(x)=\infX\supY\{f(x)+\xy-\iota_V^*(y)\}.
\end{equation*}
One wishes to invert $\infX$ and $\supY$ by means of the following
standard inf-sup theorem (see \cite{Ekeland74} for instance). We
have $\infX\supY F(x,y)=\supY\infX F(x,y)$ provided that
$\infX\supY F(x,y)\not = \pm\infty$ and
\begin{itemize}
    \item[-] $\dom F$ is a product of convex sets,
    \item[-] $x\mapsto F(x,y)$ is convex and \lsc\ for all $y,$
    \item[-] there exists $y_o$ such that $x\mapsto F(x,y_o)$ is
    coercive and
    \item[-] $y\mapsto F(x,y)$ is concave for all $x.$
\end{itemize}
Our assumptions on $f$ allow us to apply this result with
$F(x,y)=f(x)+\xy-\iota_V^*(y).$ Note that
\begin{equation}\label{A04}
    \infX f(x)>-\infty
\end{equation}
  since $f$ doesn't take the value
$-\infty$ and is assumed to be \lsc\ on a compact set. Therefore,
if $\infV f(x)<+\infty,$ we have
\begin{equation*}
    \infV
    f(x)=\supY\infX\{f(x)+\xy-\iota_V^*(y)\}=-\infY\{f^*(y)+\iota_V^*(-y)\}.
\end{equation*}

 $\bullet$ \textit{The case where $V\cap \cl\dom f=\emptyset.$}
As $\cl\dom f$ is assumed to be compact, by Hahn-Banach theorem
$\cl\dom f$ and $V$ are strictly separated: there exists $y_o\in
Y$ such that $\iota_V^*(y_o)=\sup_{x\in V} \langle x,y_o\rangle <
\inf_{\cl\dom
    f}\langle x,y_o\rangle\leq \inf_{x\in\dom f}\langle
    x,y_o\rangle.$ Hence,
\begin{equation}\label{A05}
\inf_{x\in\dom f}\{\langle x,y_o\rangle-\iota_V^*(y_o)\}>0
\end{equation}
and
\begin{eqnarray*}
 -\infY (f^*(y)+\iota^*_V(-y))
 &=& \supY\infX  \{f(x)+\xy-\iota_V(y)\}\\
   &=& \supY\inf_{x\in\dom f} \{f(x)+\xy-\iota_V(y)\} \\
   &\geq& \infX f(x) +\sup_{a>0}\inf_{x\in\dom f}\{\langle x,ay_o\rangle -\iota_V^*(ay_o)\} \\
   &=& \infX f(x) +\sup_{a>0}a\inf_{x\in\dom f}\{\langle x,y_o\rangle -\iota_V^*(y_o)\} \\
   &=& +\infty
\end{eqnarray*}
where the last equality follows from (\ref{A04}) and (\ref{A05}).
This proves that (\ref{A03}) holds with $W=V.$

$\bullet$ Finally, if (\ref{A01}) isn't satisfied, taking $W$ such
that $W\subset\inter V$ insures the strict separation of $W$ and
$\cl\dom f$ as above.
\endproof

\begin{lemma}\label{res-43}
Let the $\sigma(X,Y)$-closed convex neighbourhood $V$ of the
origin be defined by
\begin{equation}\label{eq-42}
    V=\{x\in X; \langle y_i,x\rangle\leq 1, 1\leq i\leq n\}
\end{equation}
with $n\geq 1$ and $y_1,\dots, y_n\in Y.$ Its support function
$\iota_V^*$ is $[0,\infty]$-valued, coercive and its domain is the
finite dimensional convex cone spanned by $\{y_1,\dots,y_n\}.$
More precisely, its level sets are $\{\iota_V^*\leq b\}=b\, \cv
\{y_1,\dots,y_n\}$ for each $b\geq 0$ where $\cv
\{y_1,\dots,y_n\}$ is the convex hull of $\{y_1,\dots,y_n\}.$
\end{lemma}

\proof The closed convex set $V$ is the polar set of
$N=\{y_1,\dots,y_n\}:$ $V=N^\circ.$  Let $x_1\in V$ and $x_o\in
E:=\cap_{1\leq i\leq n}\mathrm{ker\,}y_i.$ Then, $\langle
y_i,x_1+x_o \rangle = \langle y_i,x_1\rangle\leq 1.$ Hence,
$x_1+x_o\in V.$ Considering the factor space $X/E,$ we now work
within a finite dimensional vector space whose algebraic dual
space is spanned by $\{y_1,\dots,y_n\}.$

We still denote by $X$ and $Y$ these finite dimensional spaces. We
are allowed to apply the finite dimension results which are proved
in the book \cite{RW98} by Rockafellar and Wets. In particular,
one knows that if $C$ is a closed convex set in $Y,$ then the
gauge function $\gamma_C(y):=\inf\{\lambda\geq 0; y\in\lambda C\},
y\in Y$ is the support function of its polar set $C^\circ=\{x\in
X; \langle x,y\rangle\leq 1, \forall y\in C\}.$ This means that
$\gamma_C=\iota_{C^\circ}^*$ (see \cite{RW98}, Example 11.19).

As $V=(N^{\circ\circ})^{\circ}$ and $N^{\circ\circ}$ is the closed
convex hull of $N,$ i.e. $N^{\circ\circ}=\cv(N):$ the convex hull
of $N,$ we get $V=\cv(N)^\circ$  and
$$
\iota_V^*=\gamma_{\cv(N)}.
$$
In particular, for all real $b,$ $\iota_V^*(y)\leq b
\Leftrightarrow \gamma_{\cv(N)}(y)\leq b \Leftrightarrow y\in
b\,\cv(N).$ It follows that the effective domain of $\iota_V^*$ is
the convex cone spanned by $y_1,\dots, y_n$ and $\iota_V^*$ is
coercive.
\endproof

 \proof[Proof of Theorem \ref{res-41}]
 Let $\mathcal{N}(x_o)$ denote the set of
 all the neighbourhoods of $x_o\in X.$ We want to prove that
  $  \Glim k  f_k(x_o):=\sup_{U\in\mathcal{N}(x_o)}\Lim k \inf_{x\in U}
    f_k(x)=f(x_o).$
Since $f$ is \lsc, we have $f(x_o)=\sup_{U\in\mathcal{N}(x_o)}
\inf_{x\in U} f(x),$ so that it is enough to show that for all
$U\in\mathcal{N}(x_o),$ there exists $V\in\mathcal{N}(x_o)$ such
that $V\subset U$ and
\begin{equation}\label{A06}
    \Lim k \infV  f_k(x)= \infV f(x).
\end{equation}
The topology $\sigma(X,Y)$ is such that $\mathcal{N}(x_o)$ admits
the sets
$$
V=\{x\in X; |\langle y_i,x-x_o\rangle|\leq 1, i\leq n\}
$$
as a base where  $(y_1,\dots,y_n), n\geq 1$ describes the
collection of all the finite families of  vectors in $Y.$
 By Lemma \ref{res-A1}, there exists such a $V\subset U$  which
 satisfies
 $$\infV f_k(x)=-\infY h_k(y)
 \textrm{ for all $k\geq 1$ and } \infV f(x)=-\infY h(y)$$ where we denote $h_k(y)=g_k(y)+\iota_V^*(-y)$
and $h(y)=g(y)+\iota_V^*(-y),$ $y\in Y.$

Let $Z$ denote the vector space spanned by $(y_1,\dots,y_n)$ and
$h^Z_k, h^Z$ the restrictions to $Z$ of $h_k$ and $h.$ For all
$y\in Y,$ we have
\begin{equation}\label{A08}
    \iota_V^*(-y)=-\langle x_o,y\rangle +
 \iota_{V-x_o}^*(-y)
\end{equation}
and by Lemma \ref{res-43}, the effective domain of $\iota_V^*$ is
$Z.$ Therefore, to prove (\ref{A06}) it remains to show that
\begin{equation}\label{A07}
    \Lim k\infY h_k^Z(y)=\infY h^Z(y).
\end{equation}

By assumptions (b) and (d), $(h^Z_k)$ $\Gamma$-converges and
pointwise converges to $h^Z.$ Note that this $\Gamma$-convergence
is a consequence of the lower semicontinuity of the convex
conjugate $\iota_V^*$ and Proposition 6.25 of \cite{DalMaso}.

Because of assumptions (a) and (c), $(h^Z_k)$ is also a sequence
of finite convex functions which  converges pointwise to the
finite function $h^Z.$ By (\cite{Roc70}, Theorem 10.8), $(h^Z_k)$
converges to $h^Z$ uniformly on any compact subset of $Z$ and
$h^Z$ is convex.

We now consider three cases for $x_o.$

\noindent \textit{The case where }$x_o\in\dom f.$\
 We already know that $(h^Z_k)$ $\Gamma$-converges to $h^Z.$ To
 prove (\ref{A07}), it remains to check that the sequence $(h^Z_k)$ is
 equicoercive (see \cite{DalMaso}, ??).
 \\
 For all $y\in Y,$ $g(y)-\langle x_o,y\rangle\geq
 -f(x_o)$ and (\ref{A08}) imply $h^Z(y)\geq -f(x_o)+\iota_{V-x_o}^*(-y).$
 Since, $-f(x_o)>-\infty$ and $\iota_{V-x_o}^*$ is coercive (Lemma \ref{res-43}), we obtain
 that $h^Z$ is coercive. As $(h^Z_k)$
converges to $h^Z$ uniformly on any compact subset of $Z,$ it
follows that $(h^Z_k)$ is equicoercive. This proves (\ref{A07}).

\noindent \textit{The case where }$x_o\in\cl\dom f.$\
 In this case, there exists $x_o'\in\dom f$ such that $V'=x_o'+(V-x_o)/2=\{x\in X; |\langle 2y_i,x-x_o'\rangle|\leq 1, i\leq
k\}\in\mathcal{N}(x_o')$ satisfies $x_o\in V'\subset V\subset U.$
One deduces from the previous case,  that (\ref{A07})  holds true
with $V'$ instead of $V.$

\noindent \textit{The case where }$x_o\not\in\cl\dom f.$\
 As $(h^Z_k)$ $\Gamma$-converges to $h^Z,$ by (\cite{Beer93},
 Proposition 1.3.5) we have $\limsup_{n\rightarrow \infty}\infY h_k^Z(y)\leq
 \infY h^Z(y).$ As $x_o\not\in\cl\dom f,$ for any small enough $V\in
\mathcal{N}(x_o),$
 $\infY h^Z(y)=-\inf_{x\in V}f(x)=-\infty.$ Therefore, $\Lim k\infY
 h_k^Z(y)=\infY
 h(y)=-\infty$ which is (\ref{A07}).

 This completes the proof of Theorem \ref{res-41}.
\endproof

\section{$\Gamma$-convergence of minimization problems under constraints}\label{sec-Gamma2}

As the subsequent theorem demonstrates, the notion of
$\Gamma$-convergence is well-designed for minimization problems.
Let $\seq fk$ be a $\Gamma$-converging sequence of
$(-\infty,\infty]$-valued functions on a metric space $X.$ Let us
denote its limit
$$
    \Glim k f_k=f.
$$
Let $\theta:X\to Y$ be a continuous function with values in
another metric space $Y.$ Assume that for each $k\ge1,$ $f_k$ is
coercive and also that the sequence $\seq fk$ is
\emph{equi-coercive}, i.e.\  for all $a\ge 0,$
$\bigcup_{k\ge1}\{f_k\le a\}$ is relatively compact in $X.$

\begin{theorem}\label{res-10}
Under the above assumptions, the sequence of functions $\seq\psi
k$ on $Y$ which is defined by
$$
   \psi_k(y):= \inf\{f_k(x);x\in X: \theta(x)=y\},\quad y\in Y, k\ge1
$$
$\Gamma$-converges to
$$
    \psi(y):= \inf\{f(x);x\in X: \theta(x)=y\},\quad y\in Y.
$$
In particular, for any $y^*\in Y,$   there exists a sequence
$\seq{y^*}k$ in $Y$ such that $\Lim k y^*_k=y^*$ and $\Lim k
\inf\{f_k(x); x\in X:\theta(x)=y^*_k\}=\inf\{f(x); x\in
X:\theta(x)=y^*\}\in(-\infty,\infty].$

Moreover, if $y^*$ satisfies $\inf\{f(x); x\in
X:\theta(x)=y^*\}<\infty,$ then for each $k\ge1,$ the minimization
problem
$$
    f_k(x)\rightarrow\min;\quad x\in X: \theta(x)=y^*_k
$$
admits at least a minimizer $\hat{x}_k\in X.$ Any sequence
$\seq{\hat{x}}k$ of such minimizers admits at least one limit
point and any such limit point is a solution to the minimization
problem
$$
    f(x)\rightarrow\min;\quad x\in X: \theta(x)=y^*.
$$
\end{theorem}

The proof of this result which is based on Lemmas \ref{res-11} and
\ref{res-12} below, is postponed after the proofs of these lemmas.

Let $Y$ be another metric space. We consider a $\Gamma$-convergent
sequence $\seq gk$ of $[0,\infty]$-valued functions on $X\times Y$
with
$$
    \Glim k g_k=g.
$$
Let us define for each $k\ge1$ and $y\in Y,$
$$
    \psi_k(y):=\inf_{x\in X}g_k(x,y),\quad\psi(y):=\inf_{x\in
    X}g(x,y).
$$
Assume that for each $k\ge1,$ $g_k$ is coercive and also that the
sequence $\seq gk$ is equi-coercive on $X\times Y.$

\begin{lemma}\label{res-11}
Under the above assumptions on $\seq gk$,
    $
    \Glim k \psi_k=\psi
    $
in $Y.$
\end{lemma}

    \proof
Let us fix $y^*\in Y$ and prove that $\Glim k
\psi_k(y^*)=\psi(y^*).$ Since $g_k$ is assumed to be coercive, for
every $y\in Y,$ there exists $\hat{x}_{k,y}\in X$ such that
$\psi_k(y)=g_k(\hat{x}_{k,y},y).$

\par\smallskip\noindent
\textit{Lower bound.} Let $\seq yk$ be any converging sequence in
$Y$ such that $\Lim k y_k=y^*.$ we want to show that
$$
     \Liminf k \psi_k(y_k)\ge \psi(y^*).
$$
Suppose that  $\Liminf k \psi_k(y_k)<\infty,$ since otherwise
there is nothing to prove. We denote $x^*_k=\hat{x}_{k,y_k.}$
Then,
\begin{equation*}
  \Liminf k \psi_k(y_k)
  =  \Liminf k g_k(x^*_k,y_k)\overset{(a)}=  \Lim m g_m(x^*_m,y_m)\overset{(b)}= \Lim n g_n(x^*_n,y_n)
\end{equation*}
where the index $m$ at equality (a) means that we have extracted a
subsequence such that  $\Liminf k=\Lim n.$ At equality (b), once
again a new subsequence is extracted in order that $\seq{x^*}n$
converges to some limit point $x^*:$
\begin{equation*}
    \Lim n x^*_n=x^*.
\end{equation*}
The existence of  a limit point $x^*$ is insured by our
assumptions that $\Liminf k \psi_k(y_k)<\infty$ and
$\bigcup_{k\ge1}\{g_k\le a\}$ is relatively compact for all $a\ge
0.$ Now, by filling the holes in an approriate way one can
construct a sequence $\seq{\tilde{x}}k$ which admits $\seq xn$ as
a subsequence and such that $\Lim k\tilde{x}_k=x^*.$ It follows
that
\begin{equation*}
  \Liminf k \psi_k(y_k)= \Lim n g_n(x^*_n,y_n)
  \ge \Liminf k g_k(\tilde{x}_k,y_k)
  \overset{\checkmark}\ge g(x^*,y^*)\ge \psi(y^*)
\end{equation*}
which is the desired result. At the marked inequality, we have
used our assumption that $\Glim k g_k=f.$

\par\smallskip\noindent
\textit{Recovery sequence.} Under our assumptions, the
$\Gamma$-limit $g$ is coercive on $X\times Y,$ see \cite[Thm
7.8]{DalMaso}. It follows that $g(\cdot,y^*)$ is also coercive and
that there exists $\hat{x}\in\argmin g(\cdot,y^*).$ Let
$(x_k,y_k)_{k\ge1}$ be a recovery sequence of $\seq gk$ at
$(\hat{x},y^*).$ This means that $\Lim k (x_k,y_k)=(\hat{x},y^*)$
and
    $
    \Liminf k g_k(x_k,y_k)\le g(\hat{x},y^*)=\psi(y^*).
    $
We see eventually that $$\Liminf k \psi_k(y_k)\le \Liminf k
g_k(x_k,y_k)\le\psi(y^*),$$ which is the desired estimate.
\endproof
Let us fix $y^*\in Y.$ By Lemma \ref{res-11}, there exists a
sequence $\seq {y^*}k$  such that
\begin{equation}\label{eq-14}
    \Lim k y^*_k=y^*,\quad \Lim k \psi_k(y^*_k)=\psi(y^*).
\end{equation}
Let us define
$$
    \varphi_k(x):=g_k(x,y^*_k),\quad\varphi(x):=g(x,y^*),\quad
    x\in X
$$
for all $k\ge1.$ Since $g_k$ is coercive, $\varphi_k$ is also
coercive. In particular, if $\psi(y^*)=\inf_X \varphi<\infty,$ its
minimum value $\psi_k(y^*_k)=\inf_X \varphi_k$ is finite and
therefore attained at some $\hat{x}_k\in X.$

\begin{lemma}\label{res-12}
In addition to the  assumptions of Lemma \ref{res-11}, suppose
that $\inf_X \varphi<\infty.$ For each $k,$ let $\hat{x}_k$ be a
minimizer of $\varphi_k.$ Then the sequence $\seq{\hat{x}}k$
admits limit points in $X$ and any limit point is a minimizer of
$\varphi.$
\end{lemma}
    % This implies that $\Glim k\varphi_k(\hat{x})=\varphi(\hat{x})$ at each limit point $\hat{x}$ of $\seq{\hat{x}}k.$

    \proof
We have already noticed that for each $k,$ $\varphi_k$ is coercive
so that it admits one or several minimizers. Since $\Lim
k\inf_X\varphi_k=\inf_X\varphi<\infty,$ we see that $\sup_k
\inf_X\varphi_k<\infty.$ It follows from the assumed relative
compactness of $\bigcup_{k\ge1}\{g_k\le a\}$ for all $a\ge 0,$
that $\bigcup_{k\ge1}\argmin \varphi_k$ is also relatively
compact. Therefore any sequence $\seq{\hat{x}}k$ of minimizers
$\hat{x}_k\in\argmin\varphi_k$ admits at least one limit point.
\\
As $\varphi_k(\hat{x}_k)=\psi_k(y^*_k)$, we see with \eqref{eq-14}
that
$$
    \Lim k \varphi_k(\hat{x}_k)=\inf\varphi.
$$
 On the other hand, let  $\hat{x}$ be any limit point  of
$\seq{\hat{x}}k.$ There exists a subsequence (indexed by $m$ with
an abuse of notation) such that $\Lim m \hat{x}_m=\hat{x.}$
Because of the assumed $\Gamma$-limit: $\Glim kg_k=g,$ we obtain
$$
    \varphi (\hat{x}):=g(\hat{x},y^*)
    \le \Liminf m g_m(\hat{x}_m,y^*_m)
    :=\Liminf m \varphi_m(\hat{x}_m)=\Lim k \varphi_k(\hat{x}_k)=\inf\varphi.
$$
It follows that $\hat{x}$ is a minimizer of $\varphi.$
\endproof

    \proof[Proof of Theorem \ref{res-10}]
Consider the functions
$$
g_k(x,y):=f_k(x)+\iota_{\{y=\theta(x)\}},\quad (x,y)\in X\times Y,
$$
for each $k\ge1$ and
$$
g(x,y):=f(x)+\iota_{\{y=\theta(x)\}},\quad (x,y)\in X\times Y.
$$
Because of Lemmas \ref{res-11}, \ref{res-12} and \eqref{eq-14}, to
complete the proof it is enough to show that
\begin{equation}\label{eq-15}
    \Glim k g_k=g
\end{equation}
together with the coerciveness assumptions of these lemmas.

Let us begin with the coerciveness. Since for each $k\ge1,$ $f_k$
is coercive and $\theta$ is continuous, we see that for any large
enough $a,$ $\{g_k\le a\}=\{(x,y)\in X\times Y; x\in \{f_k\le a\},
y=\theta(x)\}$ is compact, i.e.\ for each $k\ge1,$ $g_k$ is
coercive. As $\seq fk$ is assumed to be equi-coercive, its
$\Gamma$-limit $f$ is coercive and it follows by the same argument
that $g$ is also coercive. We also see that
    $\bigcup_{k\ge1}\{g_k\le a\}=\{(x,y)\in X\times Y; x\in\bigcup_{k\ge1}\{f_k\le a\}, y=\theta(x)\}$
is relatively compact, i.e. $\seq gk$ is equi-coercive.

Let us prove that \eqref{eq-15} holds true. Let $(x,y)\in X\times
Y$ be fixed.  We have to prove that:
\begin{enumerate}[(i)]
    \item For any sequence $(x_k,y_k)_{k\ge1}$ such that $\Lim k(x_k,y_k)=(x,y),$\\
    $\Liminf k f_k(x_k)+\iota_{\{y_k=\theta(x_k)\}}\ge f(x)+\iota_{\{y=\theta(x)\}}.$
    \item There exists a sequence $(\tilde{x}_k,\tilde{y}_k)_{k\ge1}$ such that $\Lim k(\tilde{x}_k,\tilde{y}_k)=(x,y),$ and\\
    $\Liminf k f_k(\tilde{x}_k)+\iota_{\{\tilde{y}_k=\theta(\tilde{x}_k)\}}\le f(x)+\iota_{\{y=\theta(x)\}}.$
\end{enumerate}
Suppose first that $y\not=\theta(x).$ Then (ii) is obvious and due
to the continuity of $\theta,$ for  any sequence
$(x_k,y_k)_{k\ge1}$ such that $\Lim k(x_k,y_k)=(x,y)$ we have that
for all large enough $k,$ $\theta(x_k)\not=y_k.$ This proves (i).
\\
Now, suppose that $y=\theta(x).$ Then (i) follows from
    $\Liminf k f_k(x_k)+\iota_{\{y_k=\theta(x_k)\}}\ge
    \Liminf k f_k(x_k)\ge f(x)=f(x)+\iota_{\{y=\theta(x)\}},$
whenever $\Lim k x_k=x.$ To prove (ii), take a recovering sequence
$\seq{\tilde{x}}k$ for $\seq fk$ at $x,$ i.e.\ $\Liminf k
f_k(\tilde{x}_k)\le f(x)$ and put
$\tilde{y}_k=\theta(\tilde{x}_k),$ for each $k\ge1.$ By the
continuity of $\theta,$ $\Lim k \tilde{y}_k=y,$ so that $\Lim k
(\tilde{x}_k,\tilde{y}_k)=(x,y).$ We also have
    $
    \Liminf k f_k(\tilde{x}_k)+\iota_{\{\tilde{y}_k=\theta(\tilde{x}_k)\}}
    = \Liminf k f_k(\tilde{x}_k)\le f(x)
    = f(x)+\iota_{\{y=\theta(x)\}},
    $
which proves (ii) and completes the proof of the theorem.
\endproof

\appendix

\section{Large deviations}

\subsection*{Large deviation principle}

We refer to the monograph by Dembo and Zeitouni \cite{DZ} for a
clear exposition of the subject. Let $X$ be a polish space
furnished with its Borel $\sigma$-field. One says that the
sequence $\seq{\gamma}n$ of probability measures on $X$ satisfies
the large deviation principle (LDP for short) with scale $n$ and
rate function $I,$ if for each Borel measurable subset $A$ of $X$
we have
\begin{equation}\label{eq-02b}
    -\inf_{x\in\inter A}I(x)\overset{(\textrm{i})}\le\Liminf n \frac1n\log\gamma_n(A)
    \le \Limsup n \frac1n\log\gamma_n(A)\overset{(\textrm{ii})}\le-\inf_{x\in\cl A}I(x)
\end{equation}
where $\inter A$ and $\cl A$ are respectively the topological
interior and closure of $A$ in $X$ and the rate function
$I:X\to[0,\infty]$ is \lsc. The inequalities (i) and (ii) are
called respectively the \emph{LD lower bound} and \emph{LD upper
bound}, where LD is an abbreviation for large deviation. The LDP
is the exact statement of what was meant in previous section when
writing
$$
\gamma_n(A)\underset{n\rightarrow\infty}\asymp
\exp\left(-n\inf_{x\in A}I(x)\right)
$$
for ``all" $A\subset X.$
\\
It is sometimes too much demanding to have the upper bound
$\Limsup n \frac1n\log\gamma_n(C)\le-\inf_{x\in C}I(x)$ for all
\emph{closed} sets $C.$ One says that we have the \emph{weak LD
upper bound} if
$$
    \Limsup n\frac1n\log\gamma_n(K)\le-\inf_{x\in K}I(x)
$$
for every \emph{compact} subset $K$ of $X.$ In case $\seq\gamma n$
satisfies the LD lower bound (for all open subsets) and the weak
LD upper bound (for all compact subsets), one says that
$\seq\gamma n$ satisfies the \emph{weak LDP}.

An important instance of large deviation principle is given by the
Sanov theorem. Consider a probability measure $R\in\PX$ on the
polish space $\XX$ and furnish $\PX$ with the narrow topology
$\sigma(\PX,\CX)$ and the corresponding Borel $\sigma$-field. Let
$Z_1,Z_2,\dots$ be a sequence of independent $\XX$-valued random
variables with common law $R,$ i.e.\ $\PP(Z_i\in B)=R(B)$ for any
Borel measurable subset $B\subset \XX$  and any $i\ge1.$ In other
words $(Z_1,\dots,Z_n)\pf\PP=R^{\otimes n}$ for all $n\ge1.$

\begin{theorem}[Sanov's theorem]\label{res-01}
Under the above assumptions, the empirical measure
$$
    L^n:=\frac1n\sum_{i=1}^n\delta_{Z_i}\in\PX
$$
satisfies the LDP \footnote{This is an abuse of definition. The
correct statement should be: the sequence $((L^n)\pf\PP)_{n\ge1}$
satisfies the LDP.} in $\PX$ with scale $n.$ Its rate function  is
$H(\cdot|R):\PX\to[0,\infty],$ the relative entropy with respect
to the reference probability measure $R.$
\end{theorem}
\noindent Here, the LDP stands in $X=\PX$ and for each $n,$
$\gamma_n=(L_n)\pf\PP\in\mathrm{P}(\PX).$ For a proof of this
result, see \cite[Thm 6.2.10]{DZ}.

Next theorem states that the continuous image of a LDP is still a
LDP with the same scale.

\begin{theorem}[Contraction principle]\label{res-08}
Let $\seq\gamma n$ satisfy the LDP in $X$ with scale $n$ and rate
function $I.$ Suppose in addition that $I$ is not only \lsc, but
that it is coercive. For any continuous function $f:X\to Y$ from
$X$ to another polish space $Y$ furnished with its Borel
$\sigma$-field,
$$
    (f\pf\gamma_n)_{n\ge1}
$$
satisfies the LDP in $Y$ with scale $n$ and the rate function
$$
    J(y)=\inf\{I(x);x: f(x)=y\},\quad y\in Y.
$$
Moreover, $J$ is also coercive.
\end{theorem}
\noindent For a proof, see \cite[Thm 4.2.1]{DZ}.

Let us look at an example of application of the contraction
principle which is in the mood of this article. Consider an
independent sequence of identically distributed random paths,
i.e.\ $(Z_1,\dots,Z_n)\pf\PP=R^{\otimes n}$ where the reference
probability measure $R$ belongs to $\PO.$ The empirical measure
$L^n$ is a $\PO$-valued random variable. Now let $f$ be the
marginal projection
$$
    f(P)=(X_0,X_1)\pf P=(P_0,P_1)\in\PX\times\PX,\quad P\in\PO.
$$
It is a continuous function. This is clear when
$\Omega=C(\ii,\XX)$ and it remains true when $\Omega=D(\ii,\XX)$
($t=0,1$ being the initial and final times, $X_0$ and $X_1$ turns
out to be Skorokhod-continuous). Using the notation of the
previous section, we see that
$$
    f(L^n)=(L^n_0,L^n_1).
$$
By Sanov's theorem, the sequence of empirical measures $L^n$
satisfies the LDP in $\PO$ with scale $n$ and rate function
$H(\cdot|R).$ Applying the contraction principle with $f$ as
above, we see that $(L^n_0,L^n_1)_{n\ge1}$ satisfies the LDP in
$\PX\times\PX$ with scale $n$ and rate function
$$
    J(\mu_0,\mu_1)=\inf\{H(P|R);P\in\PO:P_0=\mu_0,P_1=\mu_1\}\in[0,\infty],
    \quad \mu_0,\mu_1\in\PX,
$$
compare \eqref{S}.

\begin{theorem}[Laplace-Varadhan principle]\label{res-15}
Suppose that $\seq{\gamma}n$ satisfy the LDP in $X$ with a
coercive rate function $I:X\to[0,\infty],$ and let $f$ be a
continuous function on $X.$ Assume further that
$$
    \Lim M \Liminf n \frac1n\log \int_X e^{nf(x)}\1_{\{f\ge
    M\}}\,\gamma_n(dx)=-\infty.
$$
Then,
    $$
    \Lim n \frac1n \log \int_X e^{nf(x)}\,\gamma_n(dx)
    =\sup_{x\in X}\{f(x)-I(x)\}.
    $$
\end{theorem}
\noindent For a proof, see \cite[Thm 4.3.1]{DZ}.

%\bibliographystyle{alpha}
%\bibliography{bib-christian}

\begin{thebibliography}{McC95}

\bibitem[ADPZ]{ADPZ10}
S.~Adams, N.~Dirr, M.~Peletier, and J.~Zimmer.
\newblock From a large-deviation proinciple to the {W}asserstein gradient flow:
  a new micro-macro passage.
\newblock Preprint. arXiv:1004.4076.

\bibitem[Att84]{Attouch84}
H.~Attouch.
\newblock {\em Variational convergence for functions and operators}.
\newblock Pitman Advanced Publishing Program. Pitman, 1984.

\bibitem[Bee93]{Beer93}
G.~Beer.
\newblock {\em Topologies on closed and closed convex sets}, volume 268 of {\em
  Mathematics and Its Applications}.
\newblock Kluwer Academic Publishers, 1993.

\bibitem[Bil68]{Bil68}
P.~Billingsley.
\newblock {\em Convergence of Probability Measures}.
\newblock Wiley Series in Probability and Mathematical Statistics. John Wiley
  \& Sons, 1968.

\bibitem[Bre91]{Bre91}
Y.~Brenier.
\newblock Polar factorization and monotone rearrangement of vector-valued
  functions.
\newblock {\em Comm. Pure Appl. Math.}, 44:375--417, 1991.

\bibitem[BS09]{BeSc08}
M.~Beiglb{\"o}ck and W.~Schachermayer.
\newblock Duality for {B}orel measurable cost functions.
\newblock {\em Trans. Amer. Math. Soc.}, to appear, 2009.

\bibitem[CZ08]{CZ08}
K.L. Chung and J.C. Zambrini.
\newblock {\em Introduction to Random Time and Quantum Randomness}.
\newblock World Scientific, 2008.

\bibitem[DG94]{DG94}
D.A. Dawson and J.~G\"artner.
\newblock Multilevel large deviations and interacting diffusions.
\newblock {\em Probab. Theory Related Fields}, 98:423--487, 1994.

\bibitem[DZ98]{DZ}
A.~Dembo and O.~Zeitouni.
\newblock {\em Large Deviations Techniques and Applications. {S}econd edition}.
\newblock Applications of Mathematics 38. Springer Verlag, 1998.

\bibitem[Eke74]{Ekeland74}
I.~Ekeland.
\newblock {\em La th\'eorie des jeux et ses applications \`a l'\'economie
  math\'ematique}.
\newblock Presses Universitaires de France, 1974.

\bibitem[F{\"o}l88]{Foe85}
H.~F{\"o}llmer.
\newblock {\em Random fields and diffusion processes, in Ecole d'Et\'e de
  Probabilit\'es de Saint-Flour XV-XVII-1985-87}, volume 1362 of {\em Lecture
  Notes in Mathematics}.
\newblock Springer, Berlin, 1988.

\bibitem[GS]{GS10}
A.~Galichon and B.~Salanie.
\newblock Matching with trade-offs: revealed preferences over competing
  characteristics.
\newblock Preprint. \texttt{http://hal.archives-ouvertes.fr/hal-00473173/en/}.

\bibitem[Kan42]{Kanto42}
L.V. Kantorovich.
\newblock On the translocation of masses.
\newblock {\em C. R. (Dokl.) Acad. Sci. URSS}, 37:199--201, 1942.

\bibitem[Kan48]{Kanto48}
L.V. Kantorovich.
\newblock On a problem of {M}onge (in {R}ussian).
\newblock {\em Uspekhi Mat. Nauk.}, 3:225--226, 1948.

\bibitem[L{\'e}o]{Leo07b}
C.~L{\'e}onard.
\newblock A saddle-point approach to the {M}onge-{K}antorovich transport
  problem.
\newblock To appear in \textit{ESAIM-COCV}. Published online: DOI:
  10.1051/cocv/2010013.

\bibitem[Mas93]{DalMaso}
G.~Dal\ Maso.
\newblock {\em An Introduction to $\Gamma$-Convergence}.
\newblock Progress in Nonlinear Differential Equations and Their Applications
  8. Birkh\"auser, 1993.

\bibitem[McC95]{McC95}
R.~McCann.
\newblock Existence and uniqueness of monotone measure-preserving maps.
\newblock {\em Duke Math. J.}, 80:309--323, 1995.

\bibitem[Mik04]{Mikami04}
T.~Mikami.
\newblock Monge's problem with a quadratic cost by the zero-noise limit of
  $h$-path processes.
\newblock {\em Probab. Theory Relat. Fields}, 129:245--260, 2004.

\bibitem[MT06]{MT06}
T.~Mikami and M.~Thieullen.
\newblock Duality theorem for the stochastic optimal control problem.
\newblock {\em Stochastic Process. Appl.}, 116:1815--1835, 2006.

\bibitem[MT08]{MT08}
T.~Mikami and M.~Thieullen.
\newblock Optimal transportation problem by stochastic optimal control.
\newblock {\em SIAM J. Control Optim.}, 47(3):1127--1139, 2008.

\bibitem[Roc97]{Roc70}
R.T. Rockafellar.
\newblock {\em Convex Analysis}.
\newblock Princeton landmarks in mathematics. Princeton University Press,
  Princeton, N.J., 1997.
\newblock First published in the Princeton Mathematical Series in 1970.

\bibitem[RT93]{RT93}
L.~R\"uschendorf and W.~Thomsen.
\newblock Note on the {S}chr\"odinger equation and {$I$}-projections.
\newblock {\em Statist. Probab. Lett.}, 17:369--375, 1993.

\bibitem[RT98]{RT98}
L.~R\"uschendorf and W.~Thomsen.
\newblock Closedness of sum spaces and the generalized {S}chr\"odinger problem.
\newblock {\em Theory Probab. Appl.}, 42(3):483--494, 1998.

\bibitem[RW98]{RW98}
R.T. Rockafellar and R.~Wets.
\newblock {\em Variational Analysis}, volume 317 of {\em Grundlehren der
  Mathematischen Wissenschaften}.
\newblock Springer, 1998.

\bibitem[Sch32]{Sch32}
E.~Schr{\"o}dinger.
\newblock Sur la th\'eorie relativiste de l'\'electron et l'interpr\'etation de
  la m\'ecanique quantique.
\newblock {\em Ann. Inst. H. Poincar\'e}, 2:269--310, 1932.

\bibitem[Vil09]{Vill09}
C.~Villani.
\newblock {\em Optimal Transport. Old and New}, volume 338 of {\em Grundlehren
  der mathematischen Wissenschaften}.
\newblock Springer, 2009.

\end{thebibliography}

\end{document}